\documentclass{article}
\usepackage{graphicx, amsmath, amssymb, bm} 
\usepackage{xcolor}
\usepackage{algorithm,algpseudocode}
\usepackage{bbm}

\usepackage[a4paper, total={6in, 9in}]{geometry}
\graphicspath{{Figures/}}

\title{Geometric decomposition of planar vector fields with a limit cycle}
\author{Lucas J. Morales Moya\\
\textit{Friedrich Miescher Institute for Biomedical Research}}
\begin{document}
\maketitle

\begin{abstract}
Mathematical modelling is a cornerstone of computational biology. While mechanistic models might describe the interactions of interest of a system, they are often difficult to study. On the other hand, abstract models might capture key features but remain disconnected from experimental manipulation. Geometric methods have been useful in connecting both approaches, although they have only been established for specific type of systems. Phenomena of biological relevance, such as limit cycles, are still difficult to study using conventional methods. In this paper, I explore an alternative description of planar dynamical systems and I present an algorithm to compute numerically the geometric structure of planar systems with a limit cycle.
\end{abstract}

\section{Introduction}
Reformulating dynamical systems in geometric terms facilitates a deeper understanding of Physics and simplifies calculations. Notable examples relate to concepts such as energy, work and potentials.
A typical example is the use of conservation of kinetic and potential energy to solve, what could be a complicated problem using kinematics, in a small number of computations.

Biology is no exception to these approaches, with the use of potential representations to describe the dynamics of bi- or multistable systems already proposed \cite{verd2014, rand2021}. Additionally, potentials have been used to describe conceptually processes such as cell differentiation, as exemplified by the Waddington landscape \cite{waddington2014}. Under this representation, the differentiation of cells into different cell types is compared to balls rolling down different valleys. Although the Waddington landscape had remained merely a useful analogy for many years, with the advent of dimensionality reduction methods and single cell sequencing, there has been a surge on the formalisation of this idea in the later years \cite{saez2022A, rand2021, verd2014, wang2011}. The appeal of such interpretation lies in the ability of reducing the intricate mess that are gene regulatory networks into simple biological variables, with a much clearer correspondence to our mental models.

When constructing geometric models on developmental progression, cyclic processes are often ignored. However, clocks are fundamental to many developmental processes, with notable examples being the segmentation clock \cite{oates2012} or the larval development of \textit{C. elegans} \cite{meeuse2020, tsiairis2021}. On top of that, some approaches still limit themselves with potential systems, those that can be described as the gradient of a scalar function $\psi$, known as the potential. In the case of two variable systems, these can be described as 
$
\bm{\dot x_{pot}} = \mathbbm{1} \nabla \psi = \nabla \psi,
$
where $\bm{\dot x} = (\frac{dx}{dt}, \frac{dy}{dt})^T = (\dot x, \dot y)^T$ is the time derivative of the variables $x$ and $y$ and $\mathbbm{1}$ is the 2x2 identity matrix.

\begin{figure}[pt]\
\centering
\includegraphics[width=11cm]{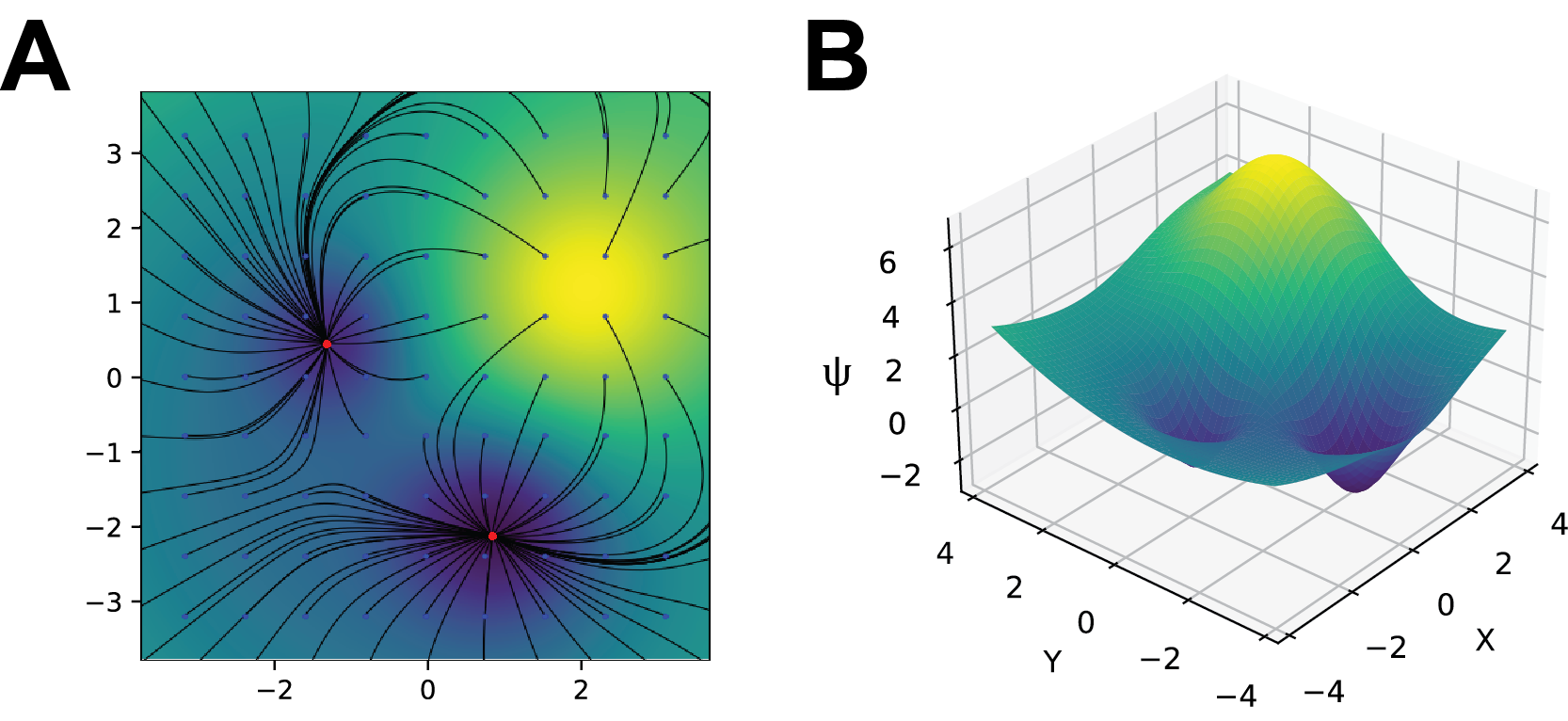} 
\caption{\textbf{Example of a potential system.}
(A) Simulation of an arbitrary potential system for different initial conditions (blue points). Each trace is shown as a black line and final simulation values are shown as red points. The value of the potential $\psi$ is shown as color, with yellow and dark blue corresponding to regions of high and low potential, respectively. Note that all trajectories converge to the regions of lowest potential, the fixed points.
(B) 3D representation of the potential of a system. Color and height represent the potential $\psi$ for different values of $x$ and $y$.} \label{pot_example_fig}
\end{figure}

Potential systems often present single point attractors, known as fixed points, to which the system evolves (see Figure \ref{pot_example_fig}). However, more complicated behaviours cannot be observed in such systems and this included notably oscillatory behaviour, as the system eventually settles at fixed points, or to infinity \cite{strogatz2018}. Rand and coworkers studied Morse-Smale potential systems and they did not pursue the analysis of phenomena such as periodic oscillations \cite{rand2021}.

However, not all systems can be described as potential, such as those exhibiting oscillatory behaviour. A different family of models are Hamiltonian systems, that can be also described as the gradient of a scalar function, the Hamiltonian ($H$), but whose vectors are rotated $\pi/2$, i.e.
$
\bm{\dot x_{rot}} = R_\perp \nabla H,
$
where $R_\perp$ is the 2x2 $\pi/2$-radians rotation matrix.

Periodic trajectories in Hamiltonian systems are neither unstable nor stable and perturbations will induce the system to follow a different orbit, in the same way that a push of a swing will increase or reduce the amplitude. These are in stark contrast with biological oscillations, which exhibits (stable) limit cycle behaviour: the system evolves to a single orbit of constant size and perturbations always return to such orbit (see Figure \ref{pot_example_fig}A) \cite{strogatz2018}. However, neither Hamiltonian nor potential systems, at least without time dependence or adding new dummy variables, can fully describe limit cycle behaviour \cite{strogatz2018, arnold2013}.

Limit cycles seem to exhibit features of both potential and Hamiltonian systems: trajectories are pulled towards a set of points, similar to potential systems, while also they orbit around fixed points, as in Hamiltonian systems. 
These properties have motivated approaches in which a vector field is decomposed in these two parts. For instance, Demongeot et al. \cite{demongeot2007I} proposed what they called the potential-Hamiltonian decomposition of a system, i.e.
$$
\bm{\dot x} = \bm{\dot x_{pot}} + \bm{\dot x_{rot}} = \mathbbm{1} \nabla \psi + R_\perp \nabla H.
$$
As noted by Suda \cite{suda2019}, this decomposition is the Helmholtz-Hodge Decomposition (HHD) in disguise, in which a vector field is divided in a curl-free and a divergence-free scalar fields, being $\psi$ and $H$ respectively.

\begin{figure}[pt]
\centering
\includegraphics[width=11cm]{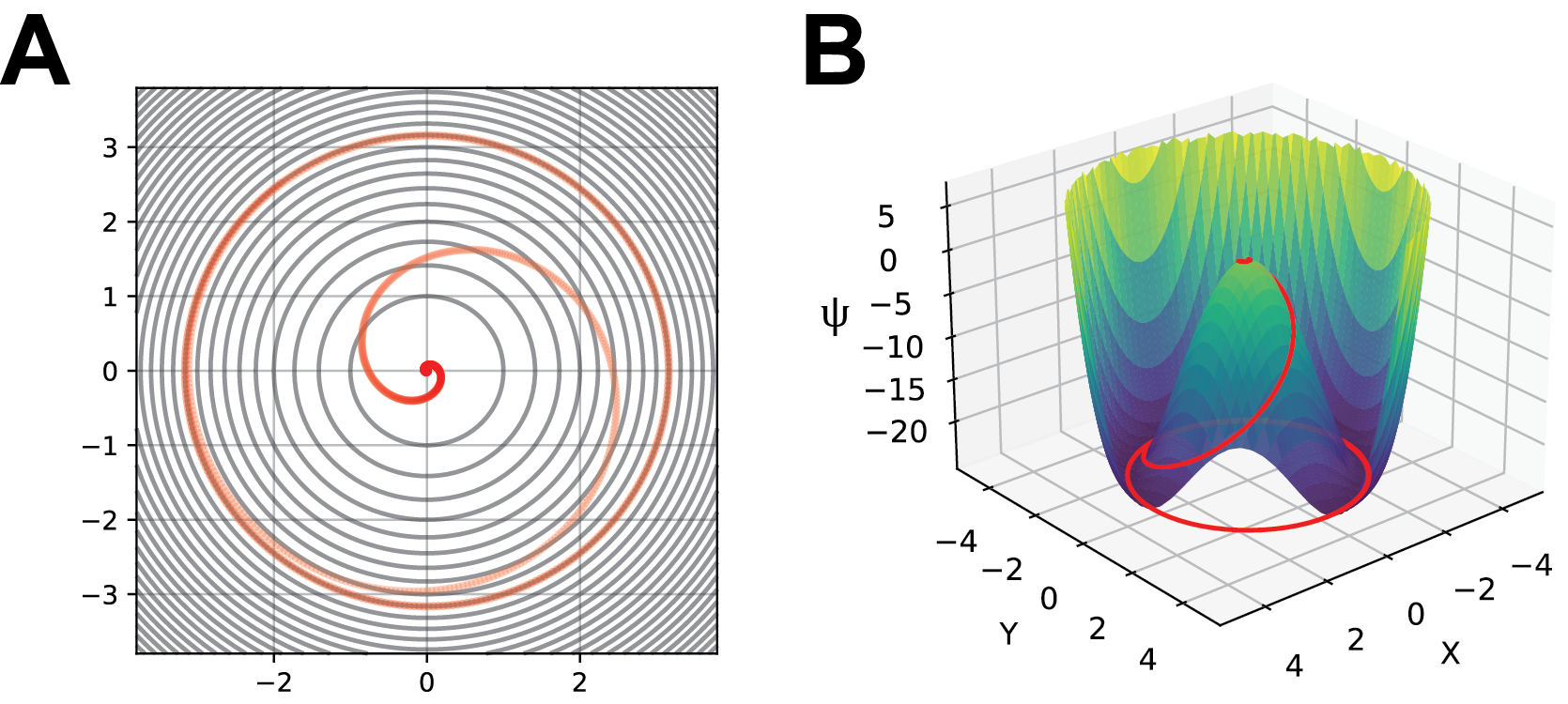} 
\caption{\textbf{The potential structure of a planar dynamical system.}
(A) The solution of Equation \ref{symm_rad} is shown in red for the initial condition $x_0 = y_0 = 3.4\cdot 10^{-3}$ (Parameter values: $\beta = 10$, $k = 0.05$, $\omega = -1$).
Energy level contours at specific points are shown as grey lines. The system converges to the trajectory given by $x^2 + y^2 = \beta$ (the circle of radius $r = \sqrt{10}$).
(B) Potential landscape for the system. Color and height represent the potential $\psi$ for different values of $x$ and $y$. Points satisfying the inequality $x + y>3$ have been removed to reveal the internal structure of the potential well. The limit cycle lies at the bottom of the potential.}
\label{symm_fig}
\end{figure}

The appeal of such decomposition lies in the identification of a potential function $\psi$ that determines how trajectories are focused upon (see Figure \ref{symm_fig}). At the surface level, the existence of a potential function highlights basins of attraction that might remain hidden in the original formulation. To an extent, this function matches the concept of the Waddington landscape and it could be seen as the formalisation of the analogy.

However, solving for both functions in the HHD requires solving two partial differential equations of second order, and thus boundary conditions are required \cite{bhatia2012}. These boundary conditions might be unclear and they might need assumptions on the behaviour of the function at different points, limiting the application of this decomposition. Even more, desired properties such as orthogonality, might not arise in this formulation \cite{suda2019, suda2020}.

Alternatively, a different decomposition was proposed in the field of Stochastic Differential Equations (SDE) \cite{ao2004}, namely the Symmetric-Antisymmetric SDE (SA-SDE) Decomposition. In this case, an SDE is decomposed in three parts: a potential function $\psi$, a dissipative symmetric matrix $S(\bm x)$, and a transversal (orthogonal) force, represented by the antisymmetric matrix $A(\bm x)$. Under this description, a SDE can be written as
$$
\left( S(\bm x) + A(\bm x) \right) \bm{\dot x} = -\nabla \psi + \xi(\bm x, t),
$$
where $\xi(\bm x, t)$ is a noise term. Under certain assumptions, one can compute the potential, as well as both matrices. However, it is worth noting that the computation of these matrices requires the definition of a diffusion matrix, based on the structure of the noise. Although some authors have claimed this could be seen as the HHD when there is no noise, this might not be always true \cite{suda2020, zhou2012}.

However, the SA-SDE decomposition can be seen as a reinterpretation of the Freidlin-Wentzell (FW) Potential for SDEs \cite{freidlin1998, wang2011, cameron2012}. In this framework, the FW potential relates to the probability of a trajectory escaping a certain basin of attraction. 

There have been other attempts at decomposing limit cycle dynamics into potential and rotational elements, these methods rely on using stochastic properties, such as diffusion \cite{wang2009, wang2010A}.
Nonetheless, this line of research is still contained within the field of SDEs, and justifications for the different assumptions are often made based on statistical arguments \cite{ao2004, wang2010B, cameron2012}. 
Still, the existence of such decompositions or a latent geometric structure should be feasible also in Ordinary Differential Equations (ODE), and should reveal the most fundamental properties of deterministic models. In other words, it should be possible to describe a deterministic dynamical system without invoking any kind of probabilistic argument or stochastic properties. Furthermore, one could analyse the system using the tools and concepts of differential geometry, such as metric tensors or differential forms \cite{frankel2011}.

Fundamentally, methods such as the HHD and the SA-SDE Decomposition rely on the gradient of the potential being zero at critical points, such as at a limit cycle or at a fixed point. However, similar results could be achieved by using Hamiltonian-like functions, which might increase monotonically and can assign a unique value to specific contours of interest.

In fact, adding more structure can help to understand the underlying geometric structure. For instance, Casimir systems are nothing but Hamiltonian systems with some non-isotropic scaling \cite{nutku1990}. A notable example is found in the Lotka-Volterra Predator-Prey Equations, where a coordinate transformation reveals a Hamiltonian structure (what Volterra called "Quantity of life" \cite{volterra1931}), i.e.
$
\bm{\dot x} = -xy R_\perp \nabla H
$, 
with $H = dx - c \ln(x) + by - a\ln (y)$.

This approach is not limited to Hamiltonian systems and it can be found in the study of potential-like systems. In particular, a general formulation of gradient descent, i.e. 
$
\bm{\dot x} = -g \nabla \psi
$, 
where $g$ is a metric tensor (a symmetric matrix that generalises the notion of distances, see \cite{misner2017} for further details) has been useful dissecting the structure of systems without orbits \cite{rand2021}.

To this day, the application of geometric decompositions to study limit cycle oscillators in biology is still in its infancy despite the large body of research on related problems and the potential to deepen our understanding of biological systems \cite{xing2010}. Additionally, the computation of solutions for these formalisms often restricted to analytical expressions of idealised systems. 

Therefore, the goal of this paper is to examine these approaches (in particular the HHD and the SA-SDED) to find novel frameworks with useful properties, as well as exploring and developing methods to compute them numerically.
Ultimately, the results developed in this paper should be the ground for future research on the geometry of biological oscillators and developmental time.

The present paper will deal with planar dynamical systems or systems of two variables, described by sufficiently smooth ODEs (at least  $C^2$ or two times differentiable) and autonomous (without explicit dependence on time), i.e.
\begin{equation}
\bm{\dot x} = f(x, y),
\end{equation}
where $\bm{\dot x} = (\frac{dx}{dt}, \frac{dy}{dt})^T = (\dot x, \dot y)^T$ is the time derivative of the variables. Additionally, I will only consider systems with a stable limit cycle (an attracting, periodic trajectory) with a single unstable fixed point at zero ($f(0,0) = 0$).

The main objective of this work is finding a reformulation of such a system in geometric terms, giving a more useful and insightful representation. In other words, we seek to find a different description of a system with two variables (mRNA, proteins, species) to a new set of variables (phase and amplitude in the case of biological oscillators) that facilitates comparison and analysis. I operate under the assumption that the problem admits such a formulation, despite not proving it. However, similar work of FW potentials can be seen as the rigorous formulation of the problem \cite{cameron2012}.

The paper is structured as follows: In Section 2, I will derive the decomposition as well as give some geometric interpretation. In Section 3, I will give further insight into the formalism using a relatively simple nonlinear model. In Section 4, I will apply the present decomposition to linear(ised) systems in an analytic way. Finally, in Section 5, I will present a novel algorithm to compute the decomposition numerically for any system with a single limit cycle enclosing a fixed point and apply to a particular case of interest.


\section{The geometric decomposition of a planar vector field} \label{derivationtion}
The proposed decomposition can be motivated and derived using the Helmholtz-Hodge Decomposition (HDD)\footnote{For the reader unfamiliar with vector calculus and its notation, I refer to the Appendix for details.}, which states that a vector field $\bm{\dot x} = f(\bm x)$ can be decomposed in two components: the gradient of a scalar potential, namely $\psi$, and a vector potential, namely $\bm{v}$, i.e.
$$
\bm{\dot x} =\nabla \psi + \nabla\times \bm{v}.
$$

In two dimensions, the vector potential $\nabla\times \bm{v}$ can be written as the gradient of another scalar function rotated by $\pi/2$, i.e. $R_\perp \nabla \phi$. Therefore, a planar vector field can be decomposed as 
\begin{equation} \label{2D_HHD}
\bm{\dot x} =\nabla \psi + R_\perp \nabla \phi,
\end{equation}
where
$$R_\perp = 
\begin{pmatrix}
0 & -1\\
1 & 0
\end{pmatrix}$$.

Generally, the rotation matrix by an angle of $\alpha$ is defined by
$$
R(\alpha) = 
\begin{pmatrix}
\cos{\alpha} & -\sin{\alpha}\\
\sin{\alpha} & \cos{\alpha}
\end{pmatrix}.
$$

\begin{figure}[pht]
\centering
\includegraphics[width = 10cm]{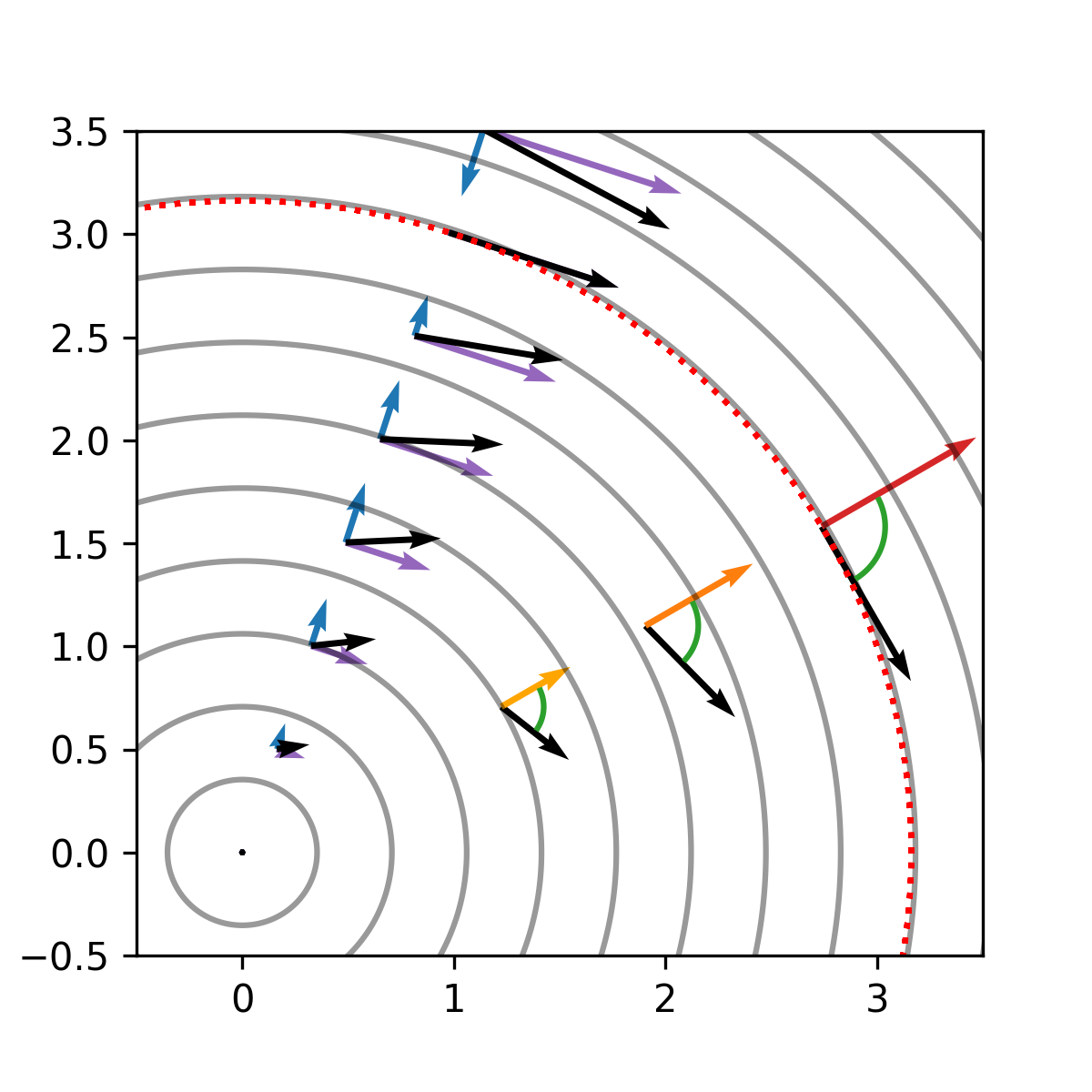} 
\caption{\textbf{Decomposition of a vector field.} Vector field (black arrows) by Equation \ref{symm_rad} (Parameter values: $\beta = 10$, $k = 0.05$, $\omega = -1$). The limit cycle is shown as a dashed red line. Some computed energy level contours for the system are shown in grey. Potential and rotational components, $\nabla \psi$ and $\nabla \phi$ respectively, are shown as blue and purple arrows. The gradient of the Hamiltonian function, $\nabla H$ is shown as warm colour arrows. Note that while $\nabla \psi$ vanishes at the limit cycle, $\nabla H$ is non-zero. The action of $\kappa$ on $\nabla H$ is shown in green. Note that at the limit cycle, it rotates $\nabla H$ by $\pi/2$.}
\label{decomp_fig}
\end{figure}

In the HHD, each one of these two new vector fields has a characteristic feature: $\nabla \psi$ is curl-free while $R_\perp \nabla \phi$ is divergence-free. In fact these properties can be shown easily (as derivatives commute for a scalar field):
$$
\nabla \times (\nabla \psi) = -\partial_y (\partial_x \psi) + \partial_x (\partial_y \psi) = 0
$$
$$
\nabla \cdot (R_\perp \nabla \phi) = \partial_x (-\partial_y \phi) + \partial_y (\partial_x \phi) = 0
$$
Ideally, both vector fields should be orthogonal to each other at all points, as the ones shown in the example in Figure \ref{decomp_fig}. To avoid clumpling, only a few selected points along a line has been shown. The black arrows are a combination of the vector fields $\nabla \psi$ and $R_\perp \nabla \phi$, shown in blue and purple, respectively. In this example, both vector fields are orthogonal, and the contribution of each field can be read directly. However, there is no guarantee that such a decomposition exist \cite{suda2019, xing2010}. 

If such orthogonal decomposition exists, then one of the vector fields can be written in terms of the other, i.e.
\begin{equation} \label{complement}
\nabla \psi = S(\bm x) R_\perp \nabla \phi,
\end{equation}
where $S(\bm x)$ is a scaling or proportionality function. In the example of Figure \ref{decomp_fig}, this is equivalent to taking the blue arrows, rotating them by $\pi/2$ and scaling them such that they match the blue arrows (see Figure \ref{decomp_fig}). Therefore, Equation \ref{2D_HHD} would become
\begin{equation} 
\bm{\dot x} = (\mathbbm{1} + S R_\perp) \nabla \psi.
\end{equation}

While this formulation could be good enough, it can lead to singularities. In some cases, the vector field might have a non-zero value, while the gradient vanishes. For instance, one cannot get the purple arrow from a blue one at the limit cycle in Figure \ref{decomp_fig}, as the gradient vanishes. These singularities have been associated to the limit cycle solutions the same way that fixed points are singularities in the solutions of ODEs \cite{yuan2017}. To avoid such issues, I can redefine $\nabla \psi$ to be the product of a scalar and a vector field, i.e.
$
\nabla \psi = p \nabla H
$ 
as well as
$
w = pS
$,
which yields
\begin{equation} \label{HHD_like_temp}
\bm{\dot x} = p\nabla H + w R_\perp \nabla H = (p \mathbbm{1} + w R_\perp) \nabla H.
\end{equation}

This equation can be written more compactly in matrix form, i.e.
\begin{equation}\label{main_equation}
\bm{\dot x} = \\
\begin{pmatrix}
p & -w \\
w & p 
\end{pmatrix}
\begin{pmatrix}
\partial_x H \\
\partial_y H
\end{pmatrix}
= \kappa \nabla H.
\end{equation}

Equation \ref{main_equation} will be the focus of this paper and the proposed geometric decomposition of a planar vector field. This formulation gives an alternative description of such vector fields, shown in the right part of Figure \ref{decomp_fig}. Under this decomposition, there is an underlying energy-like function $H$ (the different contour levels in Figure \ref{decomp_fig}). However, the vector field at each point will be rotated from the gradient of $H$ (warm coloured arrows) determined by the other geometric object $\kappa$ (green curve). 

This bilinear form $\kappa$ can present two interesting cases: when it is purely symmetric ($w = 0$) and when it is purely anti-symmetric ($p = 0$).
The first one is a special case of the gradient descent on a Riemannian manifold, where $\kappa$ is the metric tensor (traditionally written as $g$) of the manifold. Furthermore, the system is purely potential when $p = 1$ everywhere, as a vector field can be described simply as the gradient of a scalar function, i.e., $\bm {\dot x} = \nabla H$. In this case, trajectories interesect each contour level orthogonally (see Figure \ref{pot_example_fig}). The second case can be seen as a Hamiltonian system, which is in a canonical coordinate system when $w = 1$ and $\kappa$ is the symplectic form (traditionally written as $J$ or $\omega$). In this case, trajectories follow the contours of $H$ and form closed orbits.

However, the form of $\kappa$ can be greatly affected by the shape of $H$. Equation \ref{main_equation} only states that a vector field can be described by these two objects, but it says nothing about how their form. For instance, defining a different function $H$ on a Hamiltonian system whose contour levels do not match the orbits of the system will transform $\kappa$ from being skew-symmetric for all phase space to fluctuate from point to point, removing the simplicity of the decomposition and effectively rendering the approach useless. 

Therefore, the challenge of such approach is to find one or several conditions that forces the decomposition to take some form, revealing features of interest of the vector field. In particular, for a system with a limit cycle, the contour levels of $H$ might represent some type of energy function, with the limit cycle being the preferred energy level (dashed red curve in Figure \ref{decomp_fig}). Consequently, all trajectories will eventually settle into this level, and if the gradient of $H$ does not vanish there (red arrow), $\kappa$ will take the form of an anti-symmetric matrix ($p = 0$), rotating the vector field $\frac{\pi}{2}$ radians (green curve).

In summary, if a planar vector field can be decomposed as in Equation \ref{main_equation} and $\kappa$ is chosen such that it is a simple as possible, with $p$ being constant at each contour level of $H$ (and becoming zero at the limit cycle), then the proposed decomposition can yield useful insight into the structure and properties of the system. To help to exemplify this statement, I will apply this idea to a simplified system in the next section.


\section{Decomposing a stereotypical dynamical system} \label{Decomposing}
Given Equation \ref{main_equation}, I will apply it to idealised systems exhibiting limit cycle behaviour found in the literature \cite{strogatz2018, wiggins2003, guckenheimer1975, demongeot2007I, ao2004, suda2019, meeuse2020}. This exercise will reveal features of interest of these models as well as support different choices for $H$ and other constraints.

Generally, these models can be written as
\begin{equation} \label{symm}
\bm{\dot x} = 
\begin{pmatrix}
k x(\beta - x^2 - y^2) - \omega y (1 - \lambda f(x,y))\\
k y(\beta - x^2 - y^2) + \omega x (1 - \lambda f(x,y))\\
\end{pmatrix},
\end{equation}
where $\beta$ is a parameter controlling the radius (or amplitude) of the limit cycle, $k$ the limit cycle attraction strength, $\omega$ the average angular frequency of the system and, $\lambda$ the degree of frequency modulation of the system. As the system can be redimensionalised to have $k = 1$, this parameter will be assumed to be 1, unless stated otherwise.

As this model is difficult to understand as stated, a more suitable form is often sought. A coordinate transformation from Cartesian coordinates $(x, y)$ to polar coordinates $(r, \theta)$ transforms the system in Equation \ref{symm} to 
\begin{equation}
\bm{\dot r} = 
\begin{pmatrix}
r(\beta - r^2)\\
\omega (1 - \lambda f(x,y)) \\
\end{pmatrix}.
\end{equation}

This change reveals the existence of a fixed point and a limit cycle. The first equation, $\dot r = k r(\beta - r^2)$ describes the evolution along the radial coordinate and it becomes zero at the values $r = 0$, the fixed point, and $r = \sqrt{\beta}$, the limit cycle (assuming that $\beta > 0$ and $\lambda^2 < 1$).


To further understand the proposed formulation (Equation \ref{main_equation}), it is useful to explore Equation \ref{symm} with radial symmetry ($\lambda = 0$) in Cartesian coordinates, i.e.
\begin{equation} \label{symm_rad}
\bm{\dot x} = 
\begin{pmatrix}
x(\beta - x^2 - y^2) - \omega y\\
y(\beta - x^2 - y^2) + \omega x\\
\end{pmatrix}.
\end{equation}

It is easy to check that this system can be decomposed as the HHD (Equation \ref{2D_HHD}) with the following pair equations
\begin{equation} \label{symm_suda}
    \begin{cases}
      \psi = \frac{r^2}{2}(\beta-\frac{r^2}{2})\\
      \phi = \omega \frac{r^2}{2}
    \end{cases}  .
\end{equation}

This result can be found both in \cite{demongeot2007I, yuan2017, suda2019} and it is not difficult to see that both functions can give rise to a purely orthogonal description, as both are sole functions of $r$. In this case, $\nabla \psi$ is always perpendicular to the contours of $r$ while $R_\perp \nabla \phi$ is always tangent to these contours. 

As the SA-SDE Decomposition can be written in the form of Equation \ref{main_equation}, i.e.,
$$
\bm{\dot x} = (D + Q) \nabla \psi = \kappa|_{SA} \nabla \psi,
$$
it is worth comparing it to the proposed approach. The description of Equation \ref{symm_rad} using the SA-SDE Decomposition can be found in \cite{yuan2017}. In this case both the symmetric and antisymmetric matrices are $D = \mathbbm{1}$ and $Q = (\beta - r^2)^{-1} R_\perp$ \cite{yuan2017}, while the potential function $\psi$ is the same as Equation \ref{symm_suda}. Consequently, Equation \ref{symm_rad} is written as
\begin{equation} \label{symm_ao}
\bm{\dot x} = 
\begin{pmatrix}
1 & -\omega (\beta - r^2)^{-1} \\
\omega (\beta - r^2)^{-1}  & 1\\
\end{pmatrix}
\begin{pmatrix}
x(\beta - r^2)\\
y(\beta - r^2)\\
\end{pmatrix}.
\end{equation}

Note that the scalar gradient ($\nabla \psi$ in the SA-SDE Decomposition) vanishes at the origin ($r=0$) and at the limit cycle ($r = \sqrt{\beta}$), where the function $\psi$ presents an extrema (see Figure \ref{symm_fig}B). While this might not be a problem at the origin, as the vector field vanishes there, it becomes singular at the limit cycle. This is reflected on the fact that a nonzero vector field is described using a zero vector field.

Therefore, a solution to avoid such singular behaviour is imposing that the gradient can only vanish at fixed points, that is, $\nabla H = 0$ if, and only if, $\bm {\dot x} = 0$. In other words, $H$ must always increase (or decrease) along $r$. Therefore, $H$ should be ideally a function of $r$. However, these condition can still be achieved by different choices of H. For example, Equation \ref{symm_rad} can be writen in the form given by Equation \ref{main_equation} for $H = \frac{1}{2}r^2 = \frac{1}{2} (x^2 + y^2)$ as
\begin{equation} \label{symm_mine}
\bm{\dot x} = 
\begin{pmatrix}
\beta - 2H & -\omega  \\
\omega   & \beta - 2H \\
\end{pmatrix}
\begin{pmatrix}
x\\
y\\
\end{pmatrix}.
\end{equation}

As expected, the limit cycle is described by $p = 0$, which turns out to be $H_{LC} = \frac{\beta}{2} = \sqrt{\beta}$. While both Equations \ref{symm_ao} and \ref{symm_mine} describes the system using the same geometric objects, the interpretation is fundamentally different. Most notably, the scalar fields $\psi$ and $H$ characterises the system distinctly, motivating the use of different symbols to identify each of them. 

In the SA-SDE Decomposition, $\psi$ acts as a potential function, in which its extrema defines both the limit cycle and the fixed point. However, this forces $\kappa$ to become singular at the limit cycle, as the vector field does not vanish there.
On the other hand, the scalar field $H$ is closer to a Hamiltonian function, which motivates the naming of generalised Hamiltonian for this object. In this case, the function changes monotonically and only presents an extremum at the origin, responsible for the fixed point. Furthermore, $\kappa$ determines the dynamics of the system at each point and the diagonal only becomes zero at the limit cycle. Consequently, $\kappa$ must be defined everywhere, one of the conditions from the previous section.

Although both cases are different, it is worth noting that both $\psi$ and $H$ functions are related by $ \nabla \psi = (\beta - 2H) \nabla H = p(H) \nabla H$, and thus one could shift between each interpretation.


As mentioned before, the monotonicity of $H$ can be achieved in various ways. For instance, if an alternative manifold is defined, i.e., $H' =\sqrt{2H} = \sqrt{x^2 + y^2}$, Equation \ref{main_equation} becomes
\begin{equation} \label{symm_mine_r}
\bm{\dot x} = 
\begin{pmatrix}
r(\beta - r^2) & -\omega r \\
\omega r  & r(\beta - r^2) \\
\end{pmatrix}
\begin{pmatrix}
x/r\\
y/r\\
\end{pmatrix}.
\end{equation}

In both cases, $H$ and $H'$ are monotonic functions of $r$ and therefore, $\kappa$ and $\kappa'$ are well-defined everywhere, with $p$ and $p'$ being zero at the limit cycle. Although, this example shows there is not true uniqueness, the shape of the contours of both functions remain the same and the only change is due to shifts in relative distance or height. In other words, the contours are invariant. Furthermore, $p$ and $p'$ are constant at each contour level as neither function depends on the angle, i.e., $p \neq f(\theta)$. 

True uniqueness could be achieved by enforcing $w$ not to change along the complementary manifold $H^\perp$. However, it might not be possible to achieve this constraint while ensure that $p$ remains constant at each $H$. A way to circumvent this problem would be by imposing that the average value for $w$ remains constant at each level $H$, i.e.,
\begin{equation} \label{level_w}
\int_{\gamma(H)} w \ dH^\perp = w_0,
\end{equation}
for a parametrisation $\gamma(H)$ of the contour $H$. However, for the purpose of the present paper, contour uniqueness will be suficient.

Finally, it is worth comparing the proposed approach to the HHD. In section \ref{derivationtion}, Equation \ref{main_equation} was motivated using a Strictly Orthogonal Helmholtz-Hodge Decomposition (Equation \ref{2D_HHD}). In particular, each component in Equation \ref{HHD_like_temp} can be associated with a potential ($\bm{\dot x_{pot}} =\mathbbm{1} \nabla \psi$) and rotational part ($\bm{\dot x_{rot}} = R_{\perp} \nabla \psi=$), such that $\bm{\dot x} = \bm{\dot x_{pot}} + \bm{\dot x_{rot}}$. 

However, these potential and rotational elements might not necessarily correspond to the divergence-free and curl-free components of the HHD, respectively. To illustrate this, I will study the interpretation of this framework with a model where the angular velocity changes on the phase ($0< \lambda < 1$), i.e.
\begin{equation} \label{symm_phase}
\bm{\dot x} = 
\begin{pmatrix}
\beta - x^2 - y^2 & -\omega (1 - \lambda x)\\
\omega (1 - \lambda x) & \beta - x^2 - y^2 \\
\end{pmatrix}
\begin{pmatrix}
x\\
y\\
\end{pmatrix},
\end{equation}

While in this case $\bm{\dot x_{pot}}$ is curl-free, $\bm{\dot x_{rot}}$ is not divergence-free, i.e.,
\begin{equation} 
\nabla \cdot (\bm{\dot x_{rot}}) = \nabla \cdot (w R_{\perp}\nabla H) = 
\partial_x (-\omega(1-\lambda x)y) = \omega \lambda y \neq 0.
\end{equation}

Therefore, while it is still an orthogonal representation of the system it is no longer a Helmholtz-Hodge Decomposition. Additionally, while the potential component of the system remains constant along each contour, the rotational part changes magnitude at each angle. The relationship between the divergence and curl operators and the proposed decomposition are given by
\begin{equation} \label{HHD_relationships}
\begin{split}
\nabla\cdot\bm{\dot x} &= p \nabla^2 H + (\partial_H p + S \partial_{H^\perp} w) \nabla H^T\nabla H
\\
\nabla\times\bm{\dot x} &= w \nabla^2 H +\partial_H w\ \nabla H^T\nabla H
\end{split}
\end{equation}

In the case that $\partial_H w = 0$, that is, there is no change of $w$ as a function of $H$, then the curl is the scaled representation of the curvature, i.e. $\nabla\times \bm{\dot x}=w \nabla^2 H$. 

As a conclusion, despite the similarities to the HHD and the SA-SDED, the proposed approach is neither a true potential representation of the system, nor a Helmoltz-Hodge Decomposition. In fact, the proposed approach describes a planar system with a single fixed point and a limit cycle in terms of an energy-like function $H$, whose gradient vanish only at the fixed point, and a bilinear form $\kappa$, which becomes anti-symmetric ($p = 0$) at the limit cycle.

\section{Analytical decomposition of a linear system} \label{linear}
In the previous section, I have exploited the radial symmetry of the examples to illustrate how the proposed decomposition works and its relationship with other related approaches. Due to this symmetry, it was intuitive to choose $H$ to be some function of $r$. However, most of the interesting systems do no exhibit such properties and are more challenging to decompose. 

Before moving to the general case, and to fully understand the geometric meaning of Equation \ref{main_equation}, I will apply it to a linear system or, more generally, to a linearised system around a fixed point. The following results can be seen as a reformulation of those in \cite{cameron2012}, although from a purely deterministic perspective. However, final results might differ as a consequence of different constrains.

Without loss of generality, I will assume such fixed point is centered at the origin and it is described as
\begin{equation}
\bm{\dot x}|_{\bm{x} = \bm{0}} = J_0 \bm{x},
\end{equation}
where $J_0 \equiv J(\bm{x} = \bm{0})$ is the Jacobian of the system at the fixed point, defined as
\begin{equation} \label{jacobian}
J = D(\bm{\dot x}) =
\begin{pmatrix}
\partial_x \dot x & \partial_y \dot x\\
\partial_x \dot y & \partial_y \dot y\\
\end{pmatrix},
\end{equation}
where $D(\bm{v})$ is the total derivative of the vector $\bm{v}$.

Because the Jacobian is the derivative of $\bm{\dot x}$, then
\begin{equation} \label{main_jacobian}
J = D(\bm{\dot x}) = D(\kappa \nabla H) = \kappa D(\nabla H) + D(\kappa) \nabla H.
\end{equation}
However, differentiating $\kappa$ give rise to third order tensors whose treatment is more complicated using the traditional language of linear algebra. A much simpler representation can be achieved by the defining an auxiliary function\footnote{Technically, the function $\mathcal{C}$ is the embedding of the complex numbers $z = a +ib$ in the ring of 2x2 matrices over the reals, i.e. $\mathcal{C}: 	\mathbb{C} \rightarrow M_2(\mathbb{R})$}, 
$\mathcal{C}(\bm{v})$, that maps a vector $\bm{v}$ to a normal matrix, i.e.
\begin{equation} \label{complex_embedding}
\mathcal{C}(\bm{v}): 
\begin{pmatrix}
a & b
\end{pmatrix} \rightarrow
\begin{pmatrix}
a & -b \\
b & a 
\end{pmatrix},
\end{equation}
as well as the inverse operation, namely $\mathcal{C^*}(M): M \rightarrow \bm{v}$.

Equation \ref{complex_embedding} can be used to identify vectors with matrices and the other way around (via the inverse) and reformulate Equation \ref{main_equation}. By defining $\kappa = \mathcal{C}(\bm{\sigma})$, where $\bm{\sigma} = (p, w)^T$, and $Q = \mathcal{C}(\nabla H)$, Equation \ref{main_equation} can be written as
\begin{equation}
\bm{\dot x} = f(\bm{\sigma}, \nabla H) = \mathcal{C}(\bm{\sigma}) \nabla H = \mathcal{C}(\nabla H) \bm{\sigma},
\end{equation}
where $f(z_1, z_2)$ represents a type of complex multiplication (see Appendix C for some exposition of this connection). 
However, $\mathcal{C}$ will be defined to operate on $\bm{\sigma}$, as $\nabla H$ is the gradient of scalar function, and thus $\nabla \times \nabla H = 0$, which might not be true for $\bm{\sigma}$. 
Coming back to the Jacobian, and using these definitions, Equation \ref{main_jacobian} can be written as
$$
J = D(\bm{\dot x}) = \kappa D(\nabla H) + Q D(\bm{\sigma}) = \kappa \nabla^2 H + Q D(\bm{\sigma}),
$$
or in the coordinate system of reference, $(H, H^\perp)$ as
\begin{equation} \label{jacobian_decomp}
J = \kappa \nabla^2 H + Q \mathfrak S \Lambda Q^T,
\end{equation}
with 
$$
\Lambda = 
\begin{pmatrix}
\partial_{H} p & \partial_{H^\perp} p \\
\partial_{H} w & \partial_{H^\perp} w
\end{pmatrix};\quad
\mathfrak S = 
\begin{pmatrix}
1 & 0\\
0 & S 
\end{pmatrix}.
$$

This result can be checked by simple (although tedious) algebraic manipulations, and having that $H^\perp$ is orthogonal to $H$ (see Appendix: section \ref{J_proof}). Note that Equation \ref{HHD_relationships} can be obtained directly from this result. Equation \ref{jacobian_decomp} describes changes in the vector field (via the Jacobian) as the combination of those changes due to geometry and those due to dynamics. Furthermore, the term $Q \mathfrak S \Lambda Q^T$ captures the influence of changes in $\kappa$ in matrix form, as desired.

At the fixed point, the gradient vanishes and thus $Q = 0$. Alternatively, for a linear system, one can conceive that $\Lambda = 0$, that is, $\kappa$ does not change. Additionally, if $H$ can be written via Taylor series expansion, only the second order terms are relevant. Therefore, $H$ can be approximated as
$$
H|_{\bm{x}=0} \approx \frac{1}{2}\bm{x}^T F_0 \bm{x} =
\frac{1}{2} (A x^2 + B y^2 + 2 C x y),
$$
where $F_0$ is a quadratic form associated to $H$. Therefore, the (linearised) dynamics at fixed point are given by
\begin{equation} \label{jacobian_fixed}
J_0 = \kappa_0 F_0,
\end{equation}
where $\kappa_0$ is a normal matrix and $F_0$ is a positive semi-definite matrix. Therefore, Equation \ref{jacobian_fixed} can be seen as a polar-like decomposition of a $J_0$ as an operator, and thus offering a complementary representation to the traditional eigenvalue decomposition.

\begin{figure}[t]
\centering
\includegraphics[width=11cm]{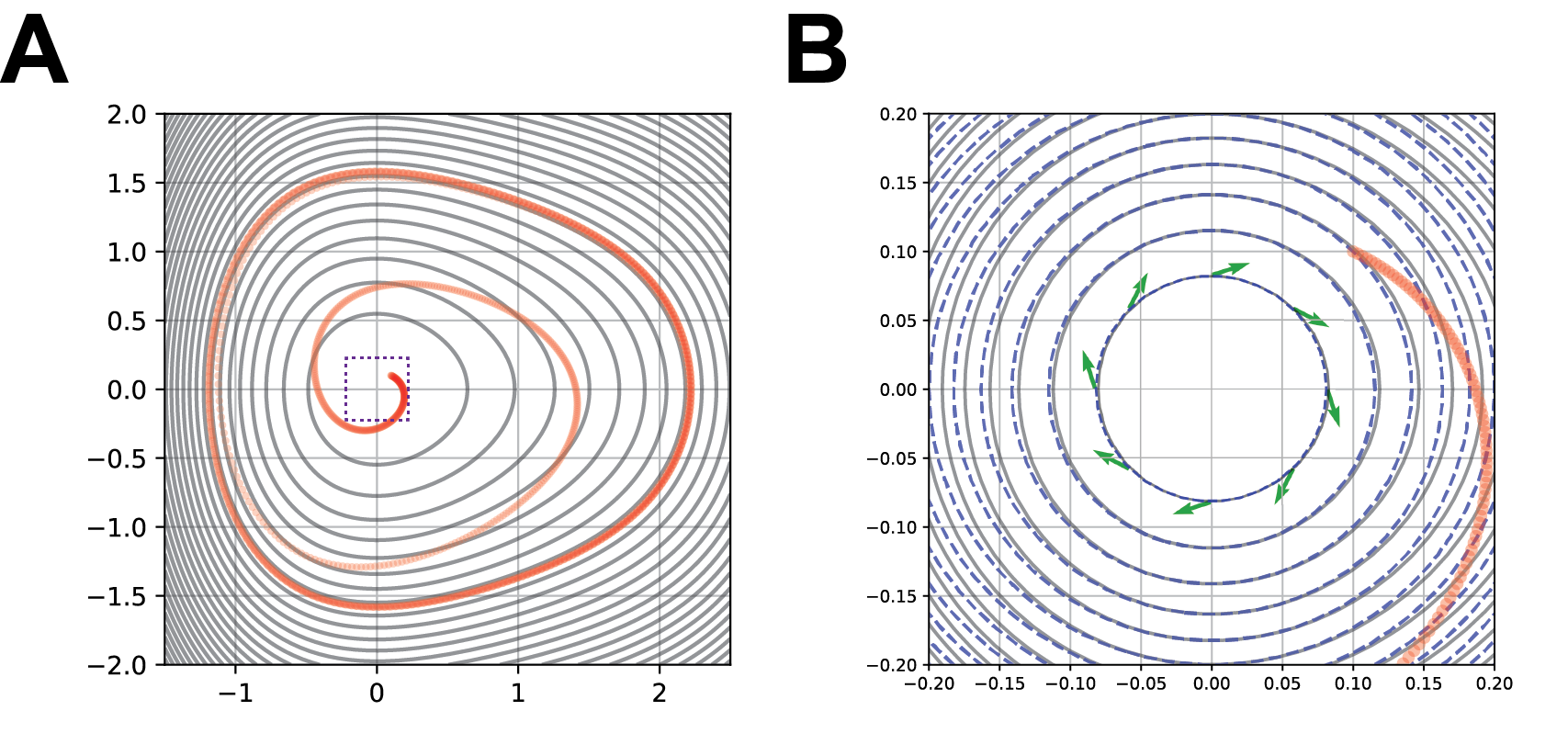} 
\caption{\textbf{Analysis and linear approximation of an asymmetric system.}
(A) The solution of Equation \ref{asymm} is shown in red for the initial condition $x_0 = y_0 = 0.1$ (Parameter values: $\beta = 10$, $k = 0.01$, $\omega = -\pi/10$, $a = 1.5$, $\mu = 2$).
Energy level contours at specific points are shown as grey lines. The system converges to the trajectory given by $H = 5$ (Equation \ref{H_asymm}). Dashed purple squared represents the region plotted in B.
(B) Linear approximation $H_0 = \frac{1}{2} (x^2 + y^2)$ (dashed blue lines) of the generalised Hamiltonian $H$ (solid grey lines) given by Equation \ref{H_asymm}. Green arrows represent the result of the action of $\kappa_0$ on $\nabla H_0$. The simulation result is shown in red.}
\label{linear_fig}
\end{figure}

A geometric interpretation of Equation \ref{jacobian_fixed} can be seen in Figure \ref{linear_fig}: $F_0$ represents the structure of the different levels near the fixed point, while $\kappa_0$ determines at which angle the vector field crosses each contour.

However, this decomposition is not unique: $J_0$ remains the same if $\kappa_0$ is divided by a factor $k$ and $F_0$ is multiplied by the same factor. This is a consequence of the covariance between both objects and it can serve as the basis of an approach based on differential geometry, in which $\kappa_0$ plays the role of a metric tensor \cite{frankel2011, chern1996, minguzzi2023}

If uniqueness is a property to be desired, then two approaches can be taken. The first one is setting the diagonal terms of the normal matrix to one and then proceed to compute the corresponding gradient term, which corresponds to the FW potential \cite{cameron2012}. However, this approach can lead to singularities for Hamiltonian systems, where $\kappa_0$ becomes undefined. Therefore, an alternative definition for this matrix is 
\begin{equation} \label{kappa_0}
\kappa_0 = \frac{1}{2}
\begin{pmatrix}
\nabla \cdot \bm{\dot x} & -\nabla \times \bm{\dot x} \\
\nabla \times \bm{\dot x} & \nabla \cdot \bm{\dot x}
\end{pmatrix}_{\bm{x}=0},
\end{equation}
This choice can be validated by some algebraic manipulations. As a consequence, $F_0 = \kappa_0^{-1} J_0$ is uniquely determined at the fixed point. Additionally, and as a result, $F_0$ gives a represents a manifold of unit mean curvature.

To illustrate the results, I choose a system with a limit cycle, with a geometry without radial symmetry whose generalised Hamiltonian is given by
\begin{equation} \label{H_asymm}
H = \frac{\mu^2}{2} (x^2 + y^2) - \frac{a}{3}x^3 + \frac{1}{4} x^4
\end{equation}
and $p = \beta - 2H$ and $w = \omega$. To ensure there is a single extremum in $H$, we require that $\mu^2 > a^2$.
In the form of Equation \ref{main_equation}, the system can be written as
\begin{equation} \label{asymm}
\dot x = 
\begin{pmatrix}
\beta - 2H & -\omega\\
\omega & \beta - 2H \\
\end{pmatrix}
\begin{pmatrix}
\mu^2 x - a x^2 + x^3\\
\mu^2 y\\
\end{pmatrix},
\end{equation}

The system has a fixed point at ${\bm{x}= \bm{0}}$. At the origin, the generalised Hamiltonian $H$ can be approximated as $H|_{\bm{x}=0} =\approx \frac{\mu^2}{2} (x^2 + y^2)$ and the Jacobian is given by
$$
J_0 = 
\begin{pmatrix}
\mu^2 \beta & -\mu^2 \omega\\
\mu^2 \omega & \mu^2 \beta  \\
\end{pmatrix}.
$$

Using Equations \ref{jacobian_fixed} and \ref{kappa_0}, one obtains the following decomposition
$$
\kappa_0 = 
\begin{pmatrix}
\mu^2 \beta & -\mu^2 \omega\\
\mu^2 \omega & \mu^2 \beta  \\
\end{pmatrix}
,\quad
F_0 = 
\begin{pmatrix}
1 & 0\\
0 & 1\\
\end{pmatrix},
$$
which gives as $H|_{\bm{x}=0} \approx \frac{1}{2} (x^2 + y^2)$, which differs by a factor of $\mu^2$ from the original estimate. However, the points belonging to a level set remains the same. In other words, there is still contour uniqueness.
If desired, one could recover the original structure by multiplying $\kappa_0$ by $\mu^{-2}$ and $F_0$ by $\mu^2$. However, both results gives the same insight into the system and thus the choice between either formulation would boil down to personal preferences. In any case, the approximation in the linear case matches the parabolic-like structure at the origin for the given $H$, getting worse as one gets farther away from the fixed point (Figure \ref{linear_fig}B).

Finally, it is worth comparing the proposed decomposition with the eigenvalue decomposition. From a traditional perspective, the fixed point is an unstable focus for positive values of $\beta$, as shown by the eigenvalues of $J_0$, i.e.
$$
\lambda_\pm = \mu^2 (\beta \pm i \omega).
$$

In this example, points $X$ belonging to a level set $H_k$, i.e. $X_k = \{ (x_j, y_j) \ | \ x_j^2 + y_j^2 = H_k \}$ will move to the next level set uniformly, at a rate defined by $\beta$. However, the points will also move alongside the level, defined by $\omega$, effectively rotating by an angle defined by $\tan^{-1}(\frac{\mu}{\beta})$. As a result, the system spirals out from the origin, as given by the eigenvalues (see Figure \ref{linear_fig}). Not only this result is recovered from the proposed decomposition, but also the expected underlying manifold $H_0$. Therefore,
the proposed decomposition is a useful tool to study linear systems.

\section{Numerical computation of the decomposition}
Up to this point, I have presented a geometric decomposition of a planar vector field with a limit cycle enclosing a fixed point and motivated its use via different examples. In particular, how it gives a description of a nonlinear system with radial symmetry as well as a linear one analytically. However, the goal of this paper goes beyond these particular cases and aims to describe the nonlinear cases with arbitrary geometries.

As giving a closed-form solution for any system might not be possible, a numerical approximation is sought. A way to get these type of solutions, is by reformulating the problem as a Hamilton-Jacobi (HJ) Equation, which can be achieved by multiplying Equation \ref{main_equation} by $(\nabla H)^T$, which gives
\begin{equation} \label{HJ}
(\nabla H)^T \dot x - p(H) (\nabla H)^T \nabla H = 0,
\end{equation}
as the $w$ term vanish due to the antisymmetry of $\kappa$.
This equation is the orthogonality condition in disguise, i.e.
\begin{equation} \label{Ortho}
(\nabla H)^T (\dot x - p(H) \nabla H) = (\nabla H)^T (w R_\perp \nabla H) = 0.
\end{equation}
This equation is a first order hyperbolic nonlinear partial differential equation and presents many challenges to solve it (for example, see Appendix C for the limitations of a naive approach). Therefore, I will present an alternative framework to get a numerical solution, exploiting the structure of the proposed approach.

To motivate this numerical approach, it is worth noticing that the contour levels of $H$ can be obtained by integrating the orthogonal gradient of $H$, i.e., $R_\perp \nabla H$, which I would call the pseudo-complement of $H$, $\nabla \hat H^\perp$. Additionally, $\nabla H$ can be compute as the product of the vector field by the inverse of $\kappa$, i.e.,
\begin{equation}
\nabla H = \kappa^{-1} \bm{\dot x},
\end{equation}
which is a direct result from Equation \ref{main_equation}. As both $\bm{\dot x}$ and $\nabla H$ do not vanish but at the fixed point, then $\kappa^{-1}$ exists for the domain of interest.

Naively, one could recover the contour levels by simply rotating all arrows by a certain angle $\alpha$ to generate a new vector field $\bm{\dot x'}$ whose vector limit cycle lies at a different contour level of $H$, i.e.,
\begin{equation} \label{rotation_alg}
\bm{\dot x'} = R(\alpha)\bm{\dot x} = R_\perp \kappa^{-1} \nabla H = \nabla \hat H^\perp.
\end{equation}
In this case, the angle $\alpha$ would indicate the relationship between $p$ and $w$ and could help to assign a value for each contour level of $H$. If this angle is sufficiently small, then the vectors belonging to the limit cycle will not be tangent to the level set, while others at a different contour might align, giving rise to a new limit cycle and revealing a new level set. By integrating the systems resulting of rotations at different angles, the whole manifold can be computed.

While this approach could work for a system such as the one given by Equation \ref{symm_rad}, this would fall short for one such as the one of Equation \ref{symm_phase}, due to the angular dependence of $w$. In other words, for this algorithm to work, $\kappa^-1$ should change according to the position in space.

As described in Section \ref{derivationtion}, the matrix $\kappa$ is an operator that scales and rotates a vector field (Equation \ref{main_equation}, Figure \ref{decomp_fig}). Consequently, Equation \ref{main_equation} can be written in terms of angles and magnitudes, i.e.,
\begin{equation} \label{set_conditions}
    \begin{cases}
      \alpha (\bm{\dot x}) = \alpha (\bm{\sigma}) + \alpha (\nabla H) \\
      |\bm{\dot x}|^2 = |\bm{\sigma}|^2 |\nabla H|^2
    \end{cases} ,
\end{equation}
where $\alpha(\bm{v}) = \tan^{-1}(v_y/v_x)$ is the angle of the vector $\bm{v}$ with respect of the $x$-axis and $|\bm{v}|^2 =(v_x^2 + v_y^2)^{-1/2}$ is the magnitude (or L2-norm) squared of the vector $\bm{v}$. The vector $\bm{\sigma} = (p, w)^T$ is the associated vector to the matrix $\kappa$, as defined in Section \ref{linear}.

The final piece of the algorithm lies in condition that $p$ is constant along a given $H$. This is sufficient for the algorithm to work. Therefore, if a given set $H_l$ is know, such as the limit cycle ($p_l = 0$), then one could compute the neighbouring sets $H_{l\pm1}$ can be computed with $p_{l\pm 1} = \pm \Delta p$.

\subsection{Level-Set Integration Algorithm}
\begin{figure}[!ht]
\centering
\includegraphics[width=12cm]{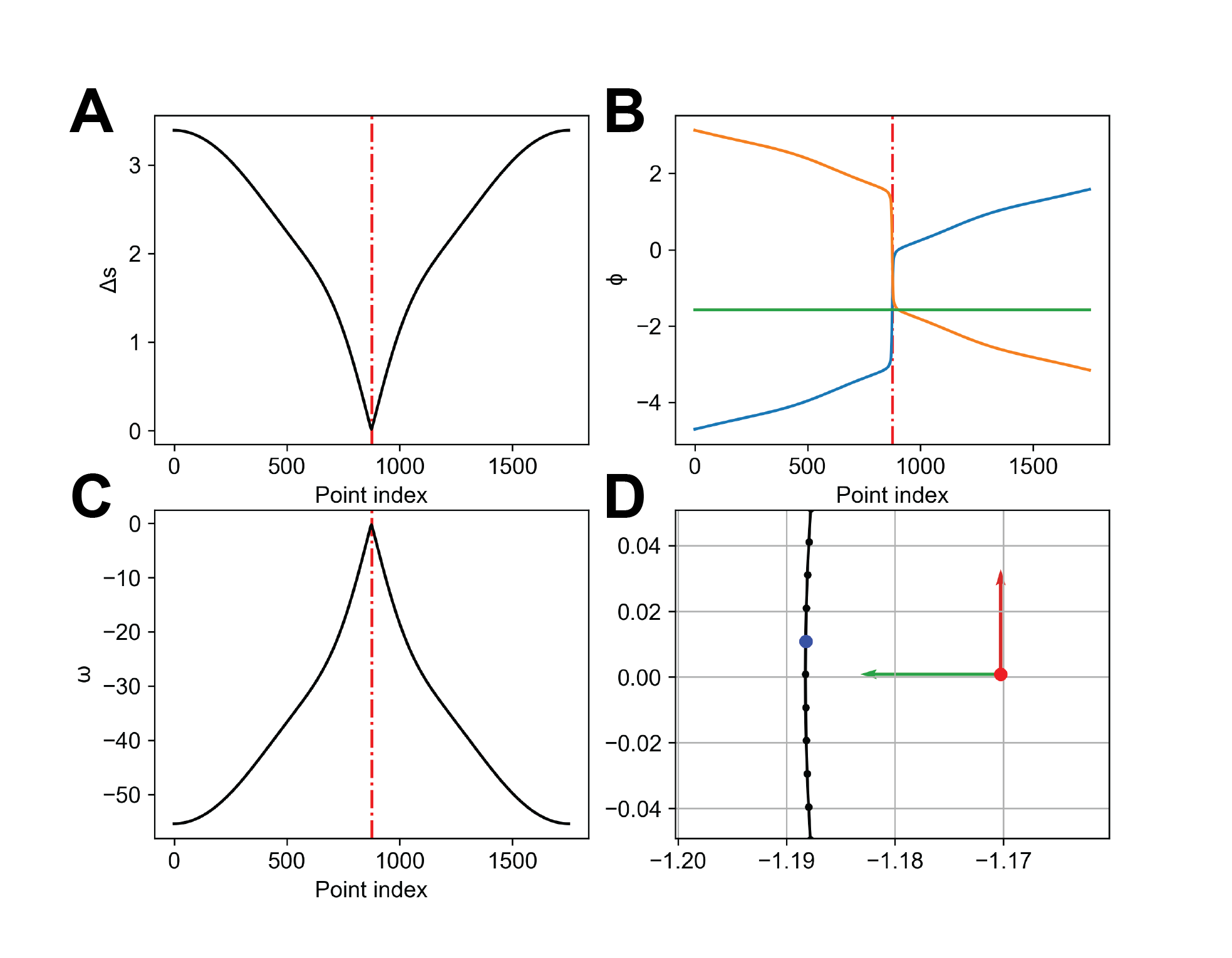} 
\caption{\textbf{Level-set integration algorithm outline.}
Computation of $\nabla H$ and $\nabla H^\perp$ for a point in relation to the limit cycle of the system of choice (Equation \ref{asymm}). Selected point of reference shown as a vertical dashed red line. Algorithm values are $p_l = 0.02$ and $\Delta H = 0.2$.
(A) Distance $\Delta s$ for the point to all the points in the reference set.
(B) Estimated angle of $\kappa$ (green line) and angular difference (blue line). Difference between both angles is shown as an orange line. The red line represent the value at which the absolute difference is minimal.
(C) Estimated values for $w$ for every point in the reference set. Computed value for $w = -0.336$ (Given value of $w_{0} = -\pi/10 \approx -0.314$).
(D) Computed vectors $\nabla \hat H$ and $\nabla \hat H^\perp$ (green and red arrows, respectively) for the selecte point (red). Reference point is shown in blue.}
\label{algorithm_outline_fig}
\end{figure}

First, a reference set $H_{l-1}$ is discretised, labelling each point as $\bm x_{l-1,j}$. A new point, outside the set and sufficiently close, namely $x_{l, k}$ is chosen. This point belongs to the level set $H_l = H_{l-1} + \Delta H$, with $\Delta H$ being chosen by the user, as well as $\Delta p = p_{l} - p_{l-1}$. Using this point, the estimated distance vector, $\Delta \hat s_{l, kj}$, from this point to the each point in the level set is given by
$$
\Delta \hat s_{l, kj} = \bm {x}_{l,k} - \bm {x}_{l-1,j} = |\Delta s_{l, kj}| \bm{v_{l, kj}},
$$
where $|\Delta s_{l, kj}|$  is the scalar distance and $\bm{v_{l, kj}}$ the normalised direction vector (see Figure \ref{algorithm_outline_fig}A).

Using the Finite Differences Method, the estimated gradient for each point, $\nabla \hat H_{l, kj}$, is computed as 
$$
\nabla \hat H_{l, kj} = \left| \frac{\Delta H}{\Delta s_{l, kj}} \right| \bm{v_{l, kj}}.
$$
It is worth noting that this gradient does not necessarily reflects any real gradient, and it is only a collection of guesses of the true value.

As both conditions from Equation \ref{set_conditions} must be satisfied, it is possible to get an estimate for $\nabla H_{l,j}$. From the second condition, $w$ can be computed as
$$
w_{l,kj} = \sqrt{|\bm {\dot x}_{l,k}|^2 |\Delta H|^{-2}  |\Delta s_{l,kj}|^2 - p_l^2}.
$$
With this result (see Figure \ref{algorithm_outline_fig}C), it is possible to compute the rotation of $\kappa_{l, kj}$, is given by
$$
\alpha_R(\kappa_{l, kj}) = \tan^{-1}\frac{w_{l, kj}}{p_l},
$$
which must be equal to the angle difference between $\bm {\dot x}_{l,k}$ and $\nabla H_{l,kj}$, given by
$$
\alpha_\Delta(\kappa_{l, kj}) = \alpha(\bm {\dot x}_{l,k}) - \alpha(\nabla H_{l,kj}).
$$

Therefore, the value for $w_{l,k}$ (and thus $\kappa_{l,k}$) is given for the one that minimizes the difference between the estimated angle of $\kappa_{l,j}$ and the angle between the vector field and the gradient, i.e.,
$$
(\alpha_\Delta(\kappa_{l, kj}) - \alpha_R(\kappa_{l, kj}))^2,
$$
as both must agree, given that the difference between both the vector field and the gradient are due to the action of $\kappa$ (see Figure \ref{algorithm_outline_fig}B). 

The estimated tangent vector to the level set of $H$, $\nabla \hat H_{l,k}$, was defined in Equation \ref{rotation_alg}. Explicitly, it is defined by
$$
\nabla \hat H_{l,k}^\perp = R_\perp \kappa_{l,k}^{-1} \bm {\dot x}_{l,k}.
$$
See Figure \ref{algorithm_outline_fig}D for an example.

To compute the rest of the points, it is necessary to use the vector $\nabla \hat H_{l,k}^\perp$ to obtain the next point in the set, ${\bm x_{l,k+1}}$. This approach is basically a Forward Euler Integration, i.e.,
\begin{equation}
\bm {x}_{l,k+1} = \bm {x}_{l,k} + \epsilon \frac{\nabla \hat H_{l,k}^\perp}{|\nabla \hat H_{l,k}^\perp|},
\end{equation}
where $\epsilon$ is a small value representing the step and chosen to be $\epsilon = 0.001$. 

Interestingly, these first points might not actually belong to $H_l$, as the choice of ${\bm x_{l,k}}$ was arbitrary. However, this algorithm is expected to eventually end up at the correct level set. Altough this result is not formally shown, and there could be cases in which this might not be true, the numerical simulations for the selected examples supported this assumption. In particular, the algorithm reached the contour level within one integrating period. To ensure the error was low enough and the contour was closed, it was integrated for three cycles and the last one was chosen to define the contour level.

\subsection{Application of the algorithm}
Before proceeding to solve complicated models with unknown geometries, I will solve the system given by Equation \ref{asymm} from the limit cycle ($p = 0$) using the Level-Set Integration Method, both inwards and outwards.

\begin{figure}[!ht]
\centering
\includegraphics[width=12cm]{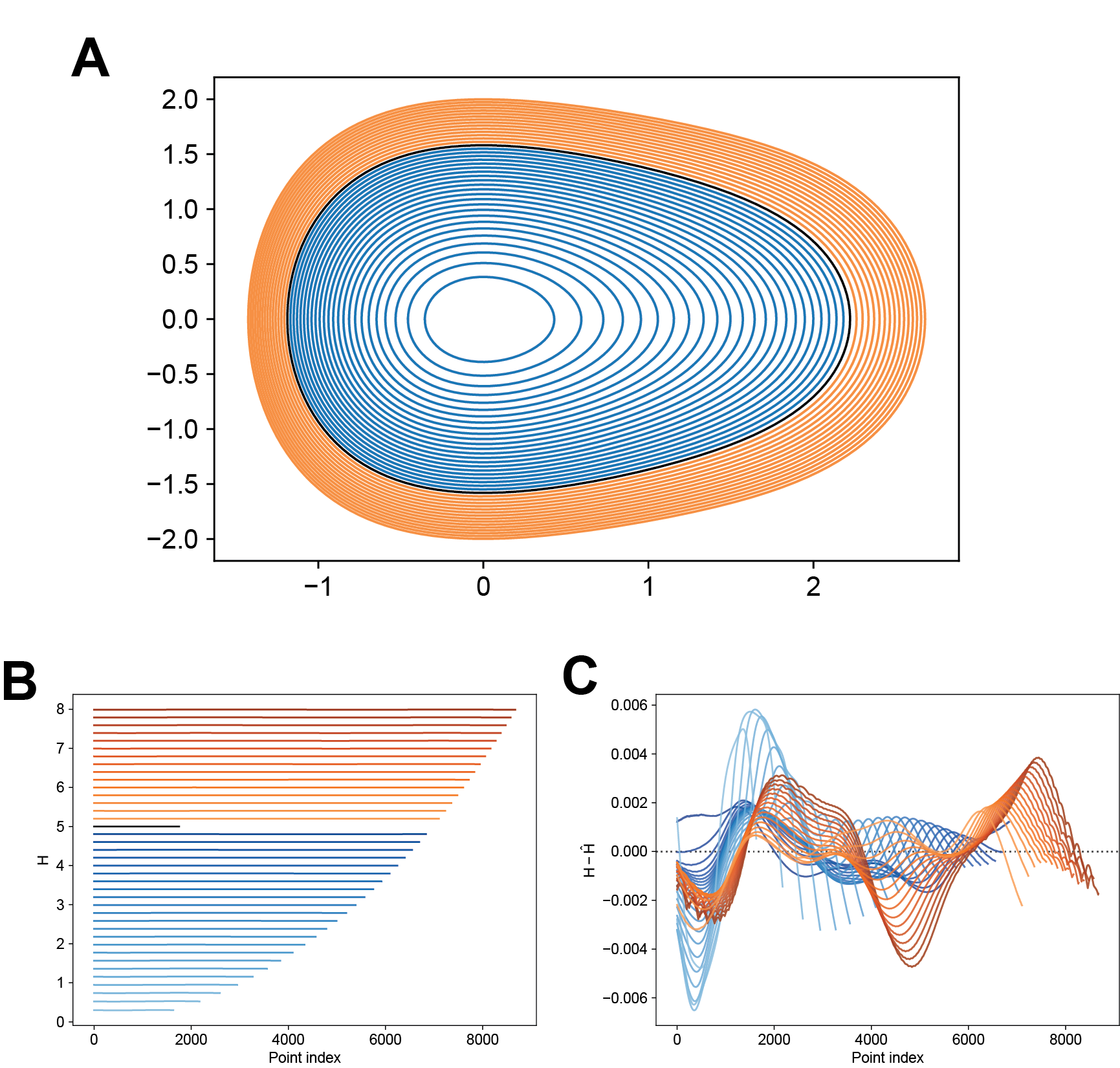} 
\caption{\textbf{Reconstruction from a limit cycle using the level-set integration method.}
(A) Reconstructed energy levels using $\Delta p = 0.02$ and $\Delta H = 0.2$ from the limit cycle (black line) of the system of choice (Equation \ref{asymm}). Inner and outer energy levels are shown as blue and orange lines, respectively.
(B) True value of $H$ for each level set.
(C) Absolute error of $H$ for each level set. Highest value of the error is given for the innermost level set (lightest blue line).}
\label{limit_cycle_fig}
\end{figure}

In both cases, the different contours levels were retrieved sucessfully for a value of $\Delta p = 0.02$ and $\Delta H = 0.2$ (see Figure \ref{limit_cycle_fig}A). In fact, the approach seems to be robustly computing the level sets from a limit cycle. As the expected value for each $H$ remains mostly constant along each contour (see Figure \ref{limit_cycle_fig}B). However, a closer looks reveal that errors might grow as the solution is propagated (see Figure \ref{limit_cycle_fig}C). 

Notably, the error increases as it gets closer to the centre, possible due to the fixed step $\Delta s$: as the manifold approaches the fixed point, $\Delta H$ decreases, and thus, a small change in $s$ lead to large changes in $H$. Therefore, an improved version of the algorithm should account for such effects to ensure that the error is bound.

\begin{figure}[!ht]
\centering
\includegraphics[width=12cm]{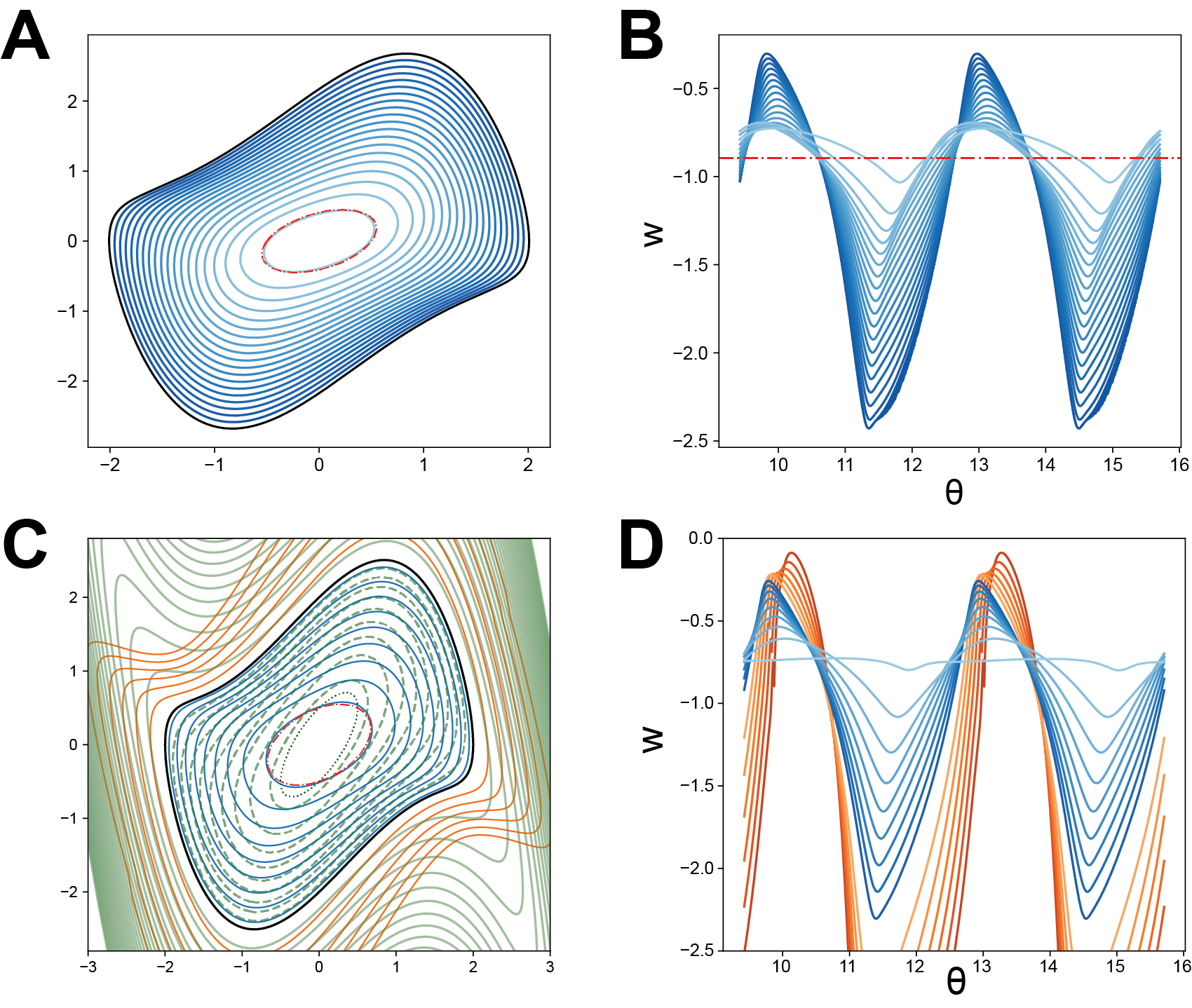} 
\caption{\textbf{Computation of the level sets for the Van der Pol oscillator and a variation.}
(A) Reconstructed energy levels using  $\Delta H = 0.1$ and $\Delta p = 0.2 \Delta H$ from the limit cycle (black line) of the system of choice (Equation \ref{vdP}). Inner level sets are shown in blue. Linear approximation for the system is shown in red.
(B) Estimated $w$ at each point for each level set. Frequency $w_0$ for the fixed point using the linearisation is shown as red dashed line.
(C) Reconstructed energy levels using  $\Delta H = 0.1$ and $\Delta p = 0.2 \Delta H$ from the limit cycle (black line) of the system of choice (Equation \ref{vdP}). Every second computed level is shown. Inner level sets are shown in blue while outer are shown in orange. Linear approximation for the system is shown in red. Computed potential from \cite{yuan2013} is shown in green. Level sets in the inner region are shown as dashed lines.
(D) Estimated $w$ at each point for each level set. Frequency $w_0$ for the fixed point using the linearisation is shown as red dashed line.}
\label{vdP_fig}
\end{figure}

A more interesting case, where the geometry is not necessarily known, is the van der Pol oscillator \cite{vanderpol1926}, one of the most prominent models exhibiting limit cycles and thus, I will use this model to test the present methodology, as well to compare it with the SA-SDED. In particular, I will use the generalisation give by \cite{yuan2013}, i.e.,
\begin{equation} \label{vdP}
\bm{\dot x} = 
\begin{pmatrix}
y\\
\mu (1-x^2) y - x + h\\
\end{pmatrix},
\end{equation}
where the traditional van der Pol oscillator is obtained for $h = 0$. The Level-Set Integration Method is able to estimate energy levels consistently, in which the limiting result matches the linear approximation at the fixed point (see Figure \ref{vdP_fig}A). As expected, this method retrieves a non-constant function for $w$ for a fixed contour level (see Figure \ref{vdP_fig}B).

To draw comparison between the approach of this paper and the one from the SA-SDE, I will use the oscillator studied by \cite{yuan2013}, with
$$
h = \frac{\mu^2}{4} x^3 - \frac{\mu^2}{16} x^5.
$$
For such system, the authors computed the gradient using the SA-SDE Decomposition as
\begin{equation} \label{yuan_potential}
\psi=\frac{u^2+v^2}{2}\left(\frac{u^2+v^2}{2} -4\right), 
\end{equation}
with $u=x$ and $v= y-\mu x+\frac{\mu}{4}x^3$.

The Jacobian of Equation \ref{vdP} is given by
\begin{equation} 
J = 
\begin{pmatrix}
0 & 1\\
-2 \mu x y - 1 + \partial_x h & \mu (1-x^2)\\
\end{pmatrix}.
\end{equation}

The gradient of the potential $\psi$ is given by
$$
\nabla\psi=\left(u^2+v^2-4\right)\ \left(
\begin{matrix}
u+\mu v\left(\frac{3}{4}x^2-1\right)\\
v\\
\end{matrix}
\right),
$$
Interestingly, a possible duality between this potential and some generalised Hamiltonian can be drawn, i.e., $\psi = H(\beta - H)$ and $\nabla \psi = (\beta - 2H) \nabla H$, in a similar fashion as done in Section \ref{Decomposing}. However, this putative generalised Hamiltonian is not necessarily the same one as the one computed using the present approach, although both algorithms agree on the limit cycle level (see Figure \ref{vdP_fig}C). Most notably, the SA-SDE Decomposition exhibits a kink at a different region than the proposed algorithm.

\begin{figure}[pht]
\centering
\includegraphics[width=11.5cm]{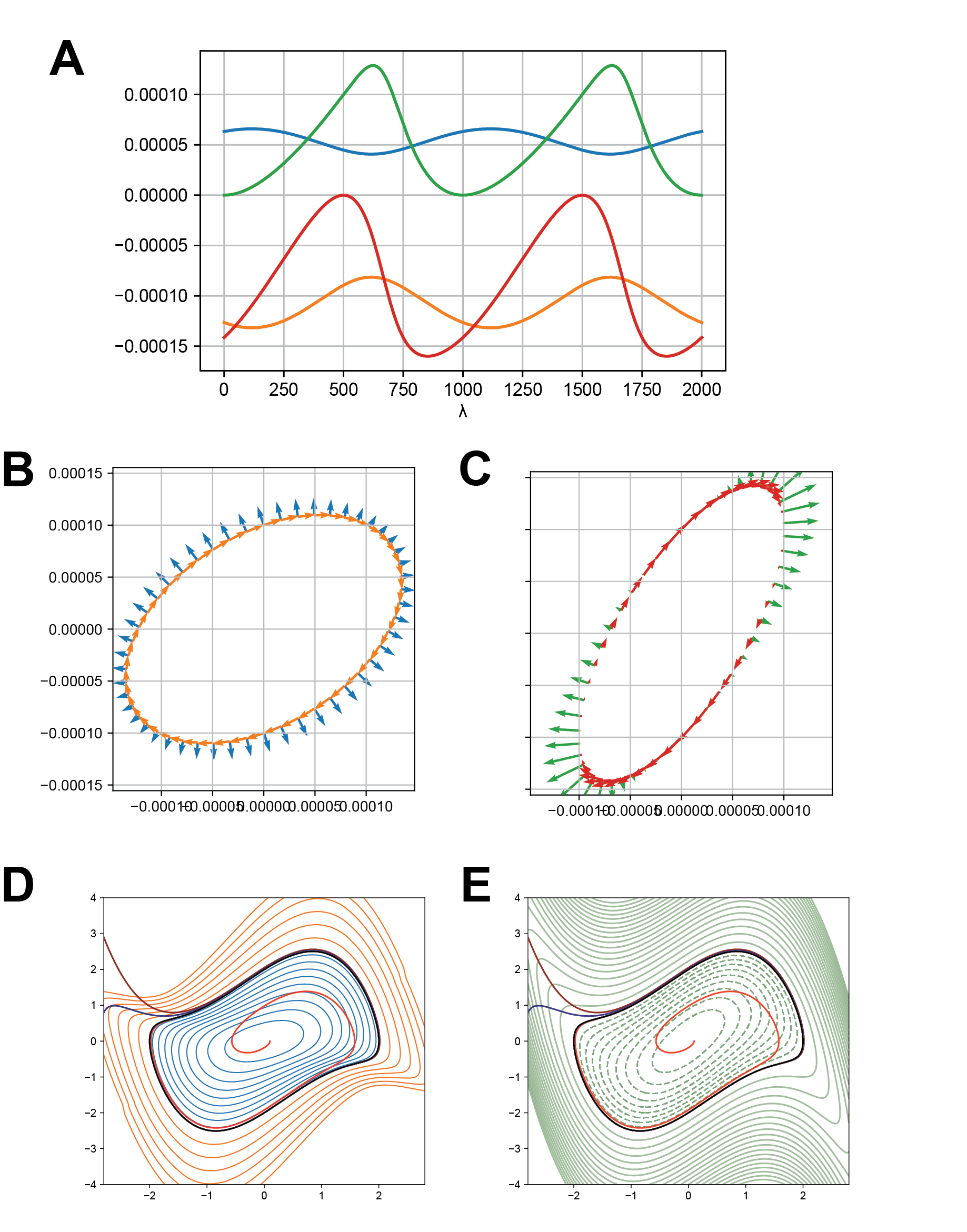} 
\caption{\textbf{Comparison between the proposed approach and the SA-SDE Decomposition.}
(A) Computed values for $p$ and $w$ along the circle for this algorithm (blue and orange, respectively) and the SA-SDE Decomposition (green and red, respectively).
(B and C) Vector field decomposition for the proposed algorithm (B) and the SA-SDE Decomposition (C).
(D and E) Comparison of three trajectories with respect to the different computed level sets of the proposed algorithm (D) and the SA-SDE Decomposition (E).}
\label{vdP_under_fig}
\end{figure}

To explore how both algorithms differ in the determination of an underlying manifold, first it is worth taking a look to the fixed point dynamics. The linearised Hamiltonian function at the fixed point can be approximated as 
$$H_0|_{R} = (1 + \mu^2) x^2 + y^2 - 2\mu xy,$$
for the approximation presented here, while the approximated linearised Hamiltonian given by the SA-SDE decomposition is 
$$H_0|_{SA} = \frac{2}{\mu^2 + 4} \left(2x^2 + (2 + \mu^2) y^2 - 2\mu xy \right)$$
These differences already reveal the angle of tilt and the eccentricity of the contour sets (see Figures \ref{vdP_fig}C, and \ref{vdP_under_fig}D and E). However, the most striking difference is the values for the potential and rotational components (see Figure \ref{vdP_under_fig}A): the proposed approach attemtps to minimise the variance along the orbit, while the SA-SDE Decomposition exhibits larger fluctuations, reaching values of zero. 

These fluctuations in the level sets are more prominent when each component is plotted spatially (see Figure \ref{vdP_under_fig}B). While the method here shows similar sized arrows along a contour, in the SA-SDE Decomposition there are points where the trajectories are tangent to the level and others where are perpendicular. This is more noticeable when a whole trajectory is plotted against the computed levels: while a trajectory starting close to the fixed point seems to cross all the level sets at the same angle in the proposed algorithm (see Figure \ref{vdP_under_fig}D), it can remain parallel to some level sets in the SA-SDE Decomposition (see Figure \ref{vdP_under_fig}E).

As expected, none of the kinks represent fully some slow manifold. However, by comparing some trajectories close to these kinks, the proposed decomposition seem to be bend towards a slow manifold chanelling trajectories into the limit cycle (see Figure \ref{vdP_under_fig}D). Nonetheless, a deeper analysis should be required to fully understand these differences, and ultimately, to distinguish between strenghts of each approach.

\section{Discussion}
First and foremost, this paper should be seen as an attempt to bridge distinct conceptual and methodological approaches from different fields to study limit cycles. Ultimately, the goal is creating a framework capable of decomposing a (planar) vector field in a set of geometric objects reporting different types of information. Although this is not a new idea and several approaches are available, such as the HHD or the SA-SDED, these are not widely used and the application to the study of Ordinary Differential Equations (ODEs) still remains largely unexplored.

In particular, the SA-SDED describes a planar system as the gradient of a potential function modified by the action of a symmetric and an antisymmetric matrices. Limit cycles arise due to extrema in such a potential function. In a nutshell, the potential function is computed using the structure of the noise in the system and relating it to the steady state distribution. Although the traditional line of reasoning for ODEs is taking the noise limit to zero, I would argue that there should be a natural way to describe a system without invoking stochasticity.

Here I have proposed a different approach to study planar differential systems and a novel algorithm to compute numerically the underlying geometry. Although the general formulation is similar to the SA-SDED \cite{ao2004} and the derivation follows a close approach to the comparison drawn by Suda between the SDE Decomposition and the Strictly Orthonormal HHD \cite{suda2019, zhou2012}, the conceptual understanding is fundamentally different, as the object of study, the scalar field, is not (necessarily) a potential, but a Hamiltonian-like function.

Both the present and the SA-SDE decompositions describe the system as a gradient system with an additional structure, a bilinear form $\kappa$. Such form actually takes the shape of the metric tensor or the symplectic form for special cases, and therefore, these decompositions can be seen as some type of unification of both approaches. As the use of both metric tensors and symplectic forms are prevalent in Physics, I would argue that the present decomposition can be of particular usefulness to study biological systems.

In Biology, there has been a wide use of potential functions to describe the dynamics of certain systems. In particular, the use of potentials to describe the temporal evolution of a biological system allows to reduce the inherent complexity of the problem to a simple geometric description. Nonetheless, the geometry of these systems is normally imposed and does not emerge naturally Mass-Action Kinetics, losing certain desired properties. 

The most prominent example of a simplified potential representation is the analogy of the Waddington landscape. Therefore, there have been efforts to formalise this concept. However, many of these approaches might overlook relevant nonlinear behaviours, such as periodic behaviour, and thus do not reflect the necessary complexity of biological dynamics. These complex behaviours beg for a more sophisticated modelling and analysis tools, which I believe this work should help to build towards.

The proposed approach defines a different set of constraints, arguing for simplicity and parsimony while avoiding the singularities at the limit cycle present in the SA-SDED. In doing so, a completely different scalar field is obtained, with a structure closer to an energy- or Hamiltonian-like function. This would allow the use of similar approaches to the ones from the symplectic geometry to be applied for problems presenting limit cycles.

Consequently, there could be new frameworks that understand biological dynamics in terms of new geometrical objects, beyond the Waddington landscape, without fully compromising the insight given by these potential representations. Furthermore, relevant nonlinear phenomena can be described using these newly developed tools.


Additionally, I have shown how such a decomposition is useful to study planar systems with a limit cycle enclosing a fixed point (centered at the origin). To build on this, I have explored linear systems using this approach and how these dynamics can be seen as the effect of an operator $\kappa$ on the gradient of a manifold, in a way akin to Polar Decompositions.

In a Polar Decomposition, any matrix can be decomposed into an isometry and a positive semi-definite matrix, which correspond to the bilinear and the quadratic matrix (or Hessian) forms, respectively. However, the bilinear form is not unitary, and thus, it is not a proper Polar Decomposition. Nonetheless, the proposed decomposition can assign a useful meaning to the positive semi-definite matrix, as it would describe the underlying geometry. Additionally, as the trace is equal to one, the average curvature is one. 

A central discussion on these methods is the existence of a unique solution. As shown for the linear case, there is uniqueness up to a constant: a combination between a scalar field and a bilinear form that satisfies Equation \ref{main_equation} for a dynamical system can be used to generate another pair of objects by scaling each one proportionally. 

In the general case, by enforcing that $p$ is just a function of $H$ and not $H^\perp$, one can achieve uniqueness up to a contour relationship. In this case, even $p$ could be described by different functions, a level set should remain the same. The difference lies in the value assigned to each level, and thus, $\partial_H p$. Although contour uniqueness should be good enough, it might be possible to enforce true uniqueness by defining \ref{level_w}. However, this has yet to be fully explored.

Alternatively, the lack of uniqueness can be seen as a consequence of the general covariance between the product of a metric tensor and a gradient \cite{frankel2011}. Under this framework, the determination of a unique solution based on constrains is just the choice of a gauge. This built-in flexibility should be seen as a feature (and possibly beneficial) and thus the uniqueness would emerge as a consequence of enforcing a perspective on the system. In this way, both the FW and the SA-SDE potentials are nothing but two valid descriptions of the same system, given complementary explanations of the system.

Interestingly, the geometry described by this metric is no longer Riemannian nor Symplectic, but Finsler \cite{chern1996, minguzzi2023}. In this case, the Finsler metric is basically a Riemannian metric without the quadratic restriction \cite{chern1996}. This opens interesting new ways to study dynamical systems, in particular biological ones, as Riemannian Geometry is still the dominant framework \cite{rand2021, verd2014}.

Finally, I have developed an algorithm to compute the proposed approach and tested it with different types of systems. In particular, I have tested a variation of the van der Pol oscillator, shown in \cite{yuan2013}, with success. In doing so, one can compare the key differences with the results from the SA-SDED and the present method.

Although a direct method could be done by solving a Hamilton-Jacobi type equation, a large number of challenges emerge in this type of problems, from unknown boundary conditions to numerical instability. Several methods have been proposed to solve these problems \cite{cameron2012, osher1988, sethian1999, yuan2013}.

The present method is a simple alternative to compute the proposed decomposition by transforming the problem into a collection of limit cycles. Although this method might not be better than other methods in the literature, this is a novel approach, which could be serve as the basis of better algorithms. This process yields sufficiently good results, although further research is required. Furthermore, the method could be modified such that integration is no longer required.

The ultimate goal of the proposed decomposition is to compute some underlying manifold that puts the system into simpler terms and reveals properties of the system. Furthermore, the present method could decompose perturbations into changes in the geometry or the dynamics and bifurcations or time-dependent effects as deformations on the geometry of the system such as ghosts \cite{strogatz2018, koch2024}

Alternatively, one can use the generalised Hamiltonian and its complement as a new coordinate system, ($H$, $H^\perp$), to describe oscillations in some natural amplitude and phase. In this case, a potential function could not be used as a coordinate system due to the vanishing gradient at the limit cycle. Putting the system in a new coordinate system resembles the traditional method of Normal Forms \cite{wiggins2003}. However, the present method does not seek to simplify the equations and the "coordinate" manifold is computed numerically and not approximated analytically.

This parametrisation of a system also could resemble isochrons \cite{winfree1980, guckenheimer1975}. In this approach, one can link points in space that (asymptotically) end up in the same point of the limit cycle, known as asymptotic phase. However, these ideas are difficult to compute and to generalise \cite{osinga2010, mauroy2013}. Therefore, the parametrisation in terms of ($H$, $H^\perp$) might offer a useful alternative to study properties of limit cycles.

In future research, it might be possible to use these tools to describe Biology in geometric terms, beyond the simple examples used in this paper. In particular, questions about biological time could be tackled in new ways using the proposed tools. For instance, the HJ equation can be written as
$$
\partial_t H(\bm{x}) = p (H(\bm{x})) |\nabla H(\bm{x})|^2,
$$
suggesting that dynamics could be solely described by the geometry of system and that the could be some notion of intrinsic time, possibly changing along the manifold. However, this is beyond the scope of this paper.

Finally, the structure and interpretation of the present decomposition, as well the ideas from this discussion, points towards concepts of Differential Geometry and ideas from Classical Mechanics and General Relativity. Therefore, using these tools could be extremely useful to understand biology. To avoid unnecessary jargon and reach a broader audience, I have avoided the language of Differential Geometry, although these results can be written using tensor algebra and calculus (For a shallow analysis, see Appendix D). Therefore, the work in the present paper should be seen as some stepping stone towards a more comprehensive and sophisticated theory of biological dynamics and time in geometric terms.

\section{Technical notes}
Simulations have been performed with Python version 3.8.5, and packages:
Numpy version 1.23.1,
Scipy version 1.10.1,
Matplotlib version 3.5.2,
Sympy version 1.12.

\section{Acknowledgements}
I would like to thank Helge Grosshans, Pavel Kos, Roman Doronin, Axel   Laborieux and Jordi García-Ojalvo for their invaluable input and suggestions. This work was supported through funding to the Grosshans lab provided in part by the Swiss National Science Foundation (Grants \#310030\_219271 and \#310030\_188487) and the Friedrich Miescher Institute for Biomedical Research through core funding through the Novartis Research Foundation.


\pagebreak
\appendix
\section{Vector calculus} \label{vector_calculus}
\subsection{Vector operators}
As the object of study are planar systems, all the results will be presented for two dimensions, although they can be extended to higher dimensions easily (with the exception of the curl operator).

The gradient of a function describes the direction and magnitude of the stepest ascent of such function. The gradient operator takes a scalar function $\psi$ and returns a vector, i.e.
$$
grad (\psi) = \nabla \psi = (\partial_x, \partial_y) \psi = (\partial_x \psi, \partial_y \psi)
$$

The divergence of a vector field is a measure of whether infinitesimally small areas expand or contract. The divergence operator takes a vector and returns a scalar, i.e.
$$
div (\bm{v}) = \nabla \cdot \bm{v} = (\partial_x, \partial_y) (v_x, v_y)^T= \partial_x v_x + \partial_y v_y
$$

The curl of a vector field is a measure of whether the vector field rotates clockwise or anticlockwise. Alternatively, it also describes path independence: if the curl is zero, the work done by the vector field between two points is the same, independently of the path taken between points. The curl operator takes a vector and returns a scalar, i.e.
$$
curl (\bm{v}) = \nabla \times \bm{v} = (-\partial_y, \partial_x) (v_x, v_y)^T= \partial_x v_y - \partial_y v_x
$$

The Laplacian of a scalar field measures the curvature of such field. It can be seen as the generalisation of the second derivative in one dimension to higher dimensions. The Laplacian computes how the gradient changes. The Laplacian operator takes a scalar function $\psi$ and returns another scalar, i.e. 
$$
lap (\psi) = \nabla^2 \psi = (\partial_x, \partial_y) (\partial_x, \partial_y)^T \psi = (\partial_x, \partial_y) (\partial_x \psi, \partial_y \psi)^T = \partial_{xx} \psi + \partial_{yy} \psi
$$

\subsection{Key identities}
The curl of the gradient is always zero, i.e.
$$
\nabla \times (\nabla \ ) = (-\partial_y, \partial_x) (\partial_x, \partial_y)^T= \partial_{xy} - \partial_{yx} = 0,
$$
as derivatives commute ($\partial_x \partial_y = \partial_y \partial_x$).

Additionally, the divergence of the curl is always zero, i.e.
$$
\nabla \cdot (\nabla \times ) = \nabla (\nabla \times ) = (\partial_x, \partial_y)(-\partial_y , \partial_x)^T = \partial_{xy} - \partial_{yx} = 0
$$
Here I have used the fact that in two dimensions, the divergence takes the form of the gradient.

\subsection{Identities in higher dimensions}
Although these are not relevant for this work, I have written them here
$$
\nabla\cdot\dot x =\nabla\cdot\left(\nabla\hat{\psi}+\nabla\times v\right)
=\nabla\cdot\nabla\hat{\psi}=\nabla^2\hat{\psi}
$$
$$
\nabla\times\dot x = \nabla\times\left(\nabla\hat{\psi}+\nabla\times v \right)=\nabla\times\nabla\times v=\nabla\cdot\nabla v -\nabla^2 v
$$

\section{Detailed proofs and demonstrations}

\subsection{Decomposition of J (Equation \ref{jacobian_decomp})} \label{J_proof}
In this section I show the decomposition of the Jacobian in its parts, i.e.
$$
J = D(\bm{\dot x}) = \kappa \nabla^2 H + Q D(\sigma) = \kappa \nabla^2 H + Q \mathfrak S \Lambda Q^T
$$

The Jacobian is given by $J = D(\bm{\dot x}) = \nabla \cdot \bm{\dot x}$. This is represented explicitly as
$$
\begin{matrix}
\begin{pmatrix}
\partial_x, \partial_y
\end{pmatrix} \\
\quad \\
\end{matrix}
\begin{pmatrix}
p\partial_x H - w\partial_y H\\
p\partial_y H + w\partial_x H
\end{pmatrix} = 
\begin{pmatrix}
\partial_x (p\partial_x H - w\partial_y H) & \partial_y (p\partial_x H - w\partial_y H)\\
\partial_x(p\partial_y H + w\partial_x H)& \partial_y(p\partial_y H + w\partial_x H)
\end{pmatrix}
$$

As $\bm{\dot x}$ is decomposed in two parts, by the chain rule, 
$$J = D(\bm{\dot x}) = D(f(\sigma, \nabla H)) = \kappa D(\nabla H) + Q D(\sigma) = J_H + J_\kappa, $$
where $J_H$ is the contribution due to the geometry and $J_\kappa$ is the contribution due to the transformation. This identity can be shown to be true by computing all the terms directly.

Changes on the geometry $J_H = \kappa \nabla^2 H$
$$
\begin{pmatrix}
p\partial_{xx} H - w\partial_{xy} H & p\partial_{yx} H - w\partial_{yy} H\\
p\partial_{xy} H + w\partial_{xx} H& p\partial_{yy} H + w\partial_{yx} H
\end{pmatrix} = 
\begin{pmatrix}
p & -w \\
w & p
\end{pmatrix}
\begin{pmatrix}
\partial_{xx} H & \partial_{yx} H \\
\partial_{xy} H & \partial_{yy} H
\end{pmatrix}
$$

Changes on the bilinear form $J_\kappa = Q D(\sigma)$
$$
\begin{pmatrix}
\partial_xp \partial_x H - \partial_xw \partial_y H & \partial_yp \partial_x H - \partial_yw \partial_y H\\
\partial_xp \partial_y H + \partial_xw \partial_x H & \partial_yp \partial_y H + \partial_yw \partial_x H
\end{pmatrix} = 
\begin{pmatrix}
\partial_{x} H & -\partial_{y} H \\
\partial_{y} H & \partial_{x} H
\end{pmatrix}
\begin{pmatrix}
\partial_{x} p & \partial_{y} p \\
\partial_{x} w & \partial_{y} w
\end{pmatrix}
$$
Ideally, I would like to express the functions $p$ and $w$ in terms of the manifold and its complement. For clarity, the complement is given by $\phi$ (instead of $H^\perp$). Therefore, the change of coordinates from $(x, y) \rightarrow (H, \phi)$ yields
$$
\begin{pmatrix}
\partial_{x} p & \partial_{y} p \\
\partial_{x} w & \partial_{y} w
\end{pmatrix}=
\begin{pmatrix}
\partial_x H \partial_H p + \partial_x \phi \partial_\phi p & 
\partial_y H \partial_H p + \partial_y \phi \partial_\phi p \\
\partial_x H \partial_H w + \partial_x \phi \partial_\phi w &
\partial_y H \partial_H w + \partial_y \phi \partial_\phi w 
\end{pmatrix},
$$
which gives $D(\sigma) = \Lambda \hat T$
$$
\begin{pmatrix}
\partial_{x} p & \partial_{y} p \\
\partial_{x} w & \partial_{y} w
\end{pmatrix}
=
\begin{pmatrix}
\partial_{H} p & \partial_{\phi} p \\
\partial_{H} w & \partial_{\phi} w
\end{pmatrix}
\begin{pmatrix}
\partial_{x} H & \partial_{y} H \\
\partial_{x} \phi & \partial_{y} \phi
\end{pmatrix}
$$
But now, I know that $\partial_x \phi = -S \partial_y H$ and $\partial_y \phi = S \partial_x H$, 
$$
\begin{pmatrix}
\partial_{x} H & \partial_{y} H \\
\partial_{x} \phi & \partial_{y} \phi
\end{pmatrix}
=\begin{pmatrix}
\partial_{x} H & \partial_{y} H \\
-S \partial_{y} H & S \partial_{x} H
\end{pmatrix} =
\begin{pmatrix}
1 & 0\\
0 & S 
\end{pmatrix}
\begin{pmatrix}
\partial_{x} H & \partial_{y} H \\
-\partial_{y} H & \partial_{x} H
\end{pmatrix},
$$
which tells that $\hat T = \mathfrak S Q^T$. Then I get
$$
J_\kappa = Q D(\sigma) = Q \Lambda \mathfrak S Q^T = Q \mathfrak S \Lambda Q^T ,
$$
with
$$
\Lambda = 
\begin{pmatrix}
\partial_{H} p & \partial_{\phi} p \\
\partial_{H} w & \partial_{\phi} w
\end{pmatrix};\quad
\mathfrak S = 
\begin{pmatrix}
1 & 0\\
0 & S 
\end{pmatrix}
$$

Without the use of the auxiliary function $\mathcal{C}$, it is more difficult to show this result using the language of linear algebra. However, using Einstein Summation Convection, Equation \ref{main_equation} can be written as
$$\dot x^i = \kappa^{ij} \partial_j H,$$
and the Jacobian is given by
$$
J_k^i = \partial_k \dot x^i = \partial_k (\kappa^{ij} \partial_j H) = \kappa^{ij} \partial_{jk} H + \partial_k\kappa^{ij} \partial_j H
$$
The change of coordinates is given by $\partial_k = \mathcal{F}_k^\beta \partial_\beta$
$$
J_k^i = \partial_k \dot x^i = \kappa^{ij} \partial_{jk} H + \mathcal{F}_k^\beta \partial_\beta \kappa^{ij} \partial_j H
$$
By redefining $\partial_\beta \kappa^{ij} = \hat \Lambda_\beta^{ij} $
$$
J_k^i = \partial_k \dot x^i = \kappa^{ij} \partial_{jk} H + \mathcal{F}_k^\beta \hat \Lambda_\beta^{ij} \partial_j H
$$

\subsection{The choice of $\kappa_0$ is correct}
I will show  
$$
J = 
\begin{pmatrix}
a & b \\
c & d
\end{pmatrix}; \quad
\kappa = 
\begin{pmatrix}
p & -w \\
w & p
\end{pmatrix}; \quad
F = 
\begin{pmatrix}
A & C \\
C & B
\end{pmatrix}
$$
Then
$$
\kappa F = \begin{pmatrix}
pA - wC & pC - wB \\
pC + wA & pB + wC
\end{pmatrix}
$$
As $J = D(\bm{\dot x})$, then
$\nabla \cdot \bm{\dot x} = a + d$ and $\nabla \times \bm{\dot x} = c - b$. Therefore,
$$
\nabla \cdot \bm{\dot x} = pA - wC + pB + wC = p(A+B)
$$$$
\nabla \times \bm{\dot x} = pC + wA - pC + wB = w(A+B) 
$$
The mean curvature is defined as $\mathfrak C = \frac{1}{2} \nabla^2 H = \frac{1}{2} (A+B)$, therefore
$$
\nabla \cdot \bm{\dot x} = 2 p \mathfrak C;  \quad 
\nabla \times \bm{\dot x} = 2 w \mathfrak C
$$

Then 
\begin{equation}
\kappa_0 = \frac{1}{2 \mathfrak{C}}
\begin{pmatrix}
\nabla \cdot \bm{\dot x} & -\nabla \cdot \bm{\dot x} \\
\nabla \cdot \bm{\dot x} & \nabla \cdot \bm{\dot x}
\end{pmatrix} \Bigg|_{\bm{x} = \bm{0}}
\end{equation}

\subsection{The decomposition can be seen as complex multiplication}
As pointed out in Section \ref{derivationtion}, the system $\bm{\dot x}$ can be thought as the rotation and scaling of the vector $\nabla H$ due to the action of $\kappa$. This is the same operation of the complex number associated to $\kappa$ on the complex number associated to $\nabla H$, resulting on the complex number associated to $\bm{\dot x}$. Explicitly,
$$
\dot x + i \dot y = (p + i w) (\partial_x H + i \partial_y H) = (p \partial_x H - w \partial_y H) + i (p\partial_y H + w \partial_x H)
$$

\subsection{Scalar fields do not transform as vector fields} \label{transformations}
If I define the transformation
$
(x, y) \rightarrow (z, y)
$, 
where $x = \mu^{-1}z$ then the system becomes
$$
\begin{pmatrix}
\dot z \\ \dot y
\end{pmatrix}=
\begin{pmatrix}
z(\beta - \mu^{-2} z^2 - y^2) - \mu \omega y\\ 
y(\beta - \mu^{-2} z^2 - y^2) + \mu^{-1}\omega  z
\end{pmatrix}
$$
which is not form-preserving, as
$$
\begin{pmatrix}
\dot z \\ \dot y
\end{pmatrix}=
\begin{pmatrix}
\beta - \mu^{-2} z^2 - y^2 & -\mu \omega\\ 
 \mu^{-1}\omega & \beta - \mu^{-2} z^2 - y^2
\end{pmatrix}
\begin{pmatrix}
z \\ y
\end{pmatrix}
$$

The transformation rule is given by
$$
\bm{\dot z} = \mathcal{F} \bm{\dot x} = \mathcal{F} \kappa_x \nabla_x H= \mathcal{F} \kappa_x \mathcal{F^T} \nabla_z H = \kappa_z \nabla_z H
$$

Is the form of $\tilde \kappa_z$ preserved in polar transformations?
$$
\mathcal{F} = 
\begin{pmatrix}
\partial_x r & \partial_y r \\
\partial_x \theta & \partial_y \theta
\end{pmatrix} 
$$

$$
\kappa_z = 
\begin{pmatrix}
\partial_x r & \partial_y r \\
\partial_x \theta & \partial_y \theta
\end{pmatrix}
\begin{pmatrix}
p & -w \\
w & p
\end{pmatrix}
\begin{pmatrix}
\partial_x r & \partial_x \theta\\
\partial_y r & \partial_y \theta
\end{pmatrix}
$$

$$
\kappa_z = 
\begin{pmatrix}
p \partial_x r + w \partial_y r & p \partial_y r - w \partial_x r\\
p \partial_x \theta + w \partial_y \theta & p \partial_y \theta - w \partial_x \theta
\end{pmatrix}
\begin{pmatrix}
\partial_x r & \partial_x \theta\\
\partial_y r & \partial_y \theta
\end{pmatrix}
$$

$$
\kappa_z = 
\begin{pmatrix}
\partial_x r (p \partial_x r + w \partial_y r) + \partial_y r (p \partial_y r - w \partial_x r)&
\partial_x \theta (p \partial_x r + w \partial_y r) + \partial_y \theta (p \partial_y r - w \partial_x r)
\\
\partial_x r(p \partial_x \theta + w \partial_y \theta) + \partial_y r (p \partial_y \theta - w \partial_x \theta) &
\partial_x \theta (p \partial_x \theta + w \partial_y \theta) + \partial_y \theta (p \partial_y \theta - w \partial_x \theta)
\end{pmatrix}
$$

$$
\kappa_z = 
\begin{pmatrix}
p (\partial_x r \partial_x r + \partial_y r \partial_y r)&
p (\partial_x \theta \partial_x r + \partial_y \theta \partial_y r)
+ w(\partial_x \theta \partial_y r - \partial_y \theta \partial_x r)
\\
p (\partial_x \theta \partial_x r + \partial_y \theta \partial_y r)
+ w(\partial_y \theta \partial_x r - \partial_x \theta \partial_y r) &
p(\partial_x \theta \partial_x \theta + \partial_y \theta \partial_y \theta)
\end{pmatrix}
$$

If we are working in polar coordinates,

$$
\mathcal{F} = 
\begin{pmatrix}
\partial_x r & \partial_y r \\
\partial_x \theta & \partial_y \theta
\end{pmatrix} = 
\begin{pmatrix}
x/r & y/r \\
-y/r^2 & x/r^2
\end{pmatrix}
$$

$$
\kappa_z = 
\begin{pmatrix}
x/r & y/r \\
-y/r^2 & x/r^2
\end{pmatrix}
\begin{pmatrix}
p & -w \\
w & p
\end{pmatrix}
\begin{pmatrix}
x/r & -y/r^2 \\
y/r & x/r^2
\end{pmatrix}
$$

$$
\kappa_z = 
\begin{pmatrix}
p\frac{x^2+ y^2}{r^2}&
p \frac{-yx + xy}{r^3}
+ w \frac{-y^2 - x^2}{r^3}
\\
p \frac{-yx + xy}{r^3}
+ w \frac{y^2 + x^2}{r^3}&
p\frac{y^2+ x^2}{r^4}
\end{pmatrix}
$$

$$
\kappa_z = 
\begin{pmatrix}
p & -w/r \\
w/r & p / r^2
\end{pmatrix}
$$
This does not seem to be form invariant.

\section{Solving a Hamilton-Jacobi equation via a naive Finite Differences Method} \label{HJ_sec}
A HJ equation can be obtained multiplying Equation \ref{main_equation} by $(\nabla H)^T$, which gives
\begin{equation} \label{HJ_Ap}
(\nabla H)^T \dot x - p(H) (\nabla H)^T \nabla H = 0,
\end{equation}
as the $w$ term vanish due to the antisymmetry of $\kappa$.
This equation is the (geometric) orthogonality condition in disguise, i.e.
\begin{equation} \label{Ortho_Ap}
(\nabla H)^T (\dot x - p(H) \nabla H) = (\nabla H)^T (w R_\perp \nabla H) = 0
\end{equation}
It is worth noting that for the two dimensional case, both elements might be non-zero. For instance, a limit cycle emerges naturally from this equation. For $p = 0$, Equation \ref{HJ} becomes
$$
\dot x \  \partial_x H + \ \dot y \ \partial_y H = 0,
$$
which gives $\partial_x H = -\dot y$ and $\partial_y H = \dot x$. Such result reveals that limit cycle solutions are perpendicular to a specific contour. The fact that $p = 0$ is a consequence of the limit cycle set orthogonality.

Equation \ref{HJ} is well-known to pose problems when integrating numerically, as naive, forward algorithms tend to be numerically unstable and might exhibit unrelated oscillations as the simulation progresses. 

Additionally, the choice of boundary conditions is not trivial: while first order methods might only require a single boundary condition, which could be our reference fixed point, higher order methods (such as Lax-Friedrich or Lax-Wendroff) might requires a second boundary condition due to the addition of artificial diffusion.

Finally, the choice of $p(H)$ must be explicit. In the case of a system with a fixed point surrounded by a limit cycle, a linear function might be sufficient, i.e. 
\begin{equation} \label{ref_pH}
p(H) = p_0 - q H
\end{equation}
This formulation is agnostic to the shape and formulation of $H$. As a consequence, it is possible to compute the gradient $\psi$ by setting $p = 1$.

\subsection{Radially symmetric system}
To illustrate the limitations of the method and how it can fail at the region close to the limit cycle ($p \rightarrow 0$), I will attempt to solve Equation \ref{symm} with $\lambda =0$ using the Hamilton-Jacobi equation (Equation \ref{HJ}) using a forward finite difference method.

The solution should be able to be propagated from our reference fixed point forward. Therefore, a polar coordinate system is more suitable. Equation \ref{HJ} becomes
\begin{equation} \label{full_polar}
\dot r \partial_r H + r^{-1} \dot \theta \partial_\theta H - p(H) [(\partial_r H)^2 + r^{-2}  (\partial_\theta H)^2] = 0
\end{equation}
For a radially symmetric system, in which $\partial_\theta H = 0$, and assuming $p(H)$ to be linear, Equation \ref{full_polar} becomes
\begin{equation} \label{radial_polar}
\dot r \partial_r H - (p_0 - q H) (\partial_r H)^2 = 0
\end{equation}
Directly, the trivial solution $\partial_r H = 0$ can be removed to simplify calculations, in this one-dimensional case. 

\begin{figure}[t]
\centering
\includegraphics[width=6cm]{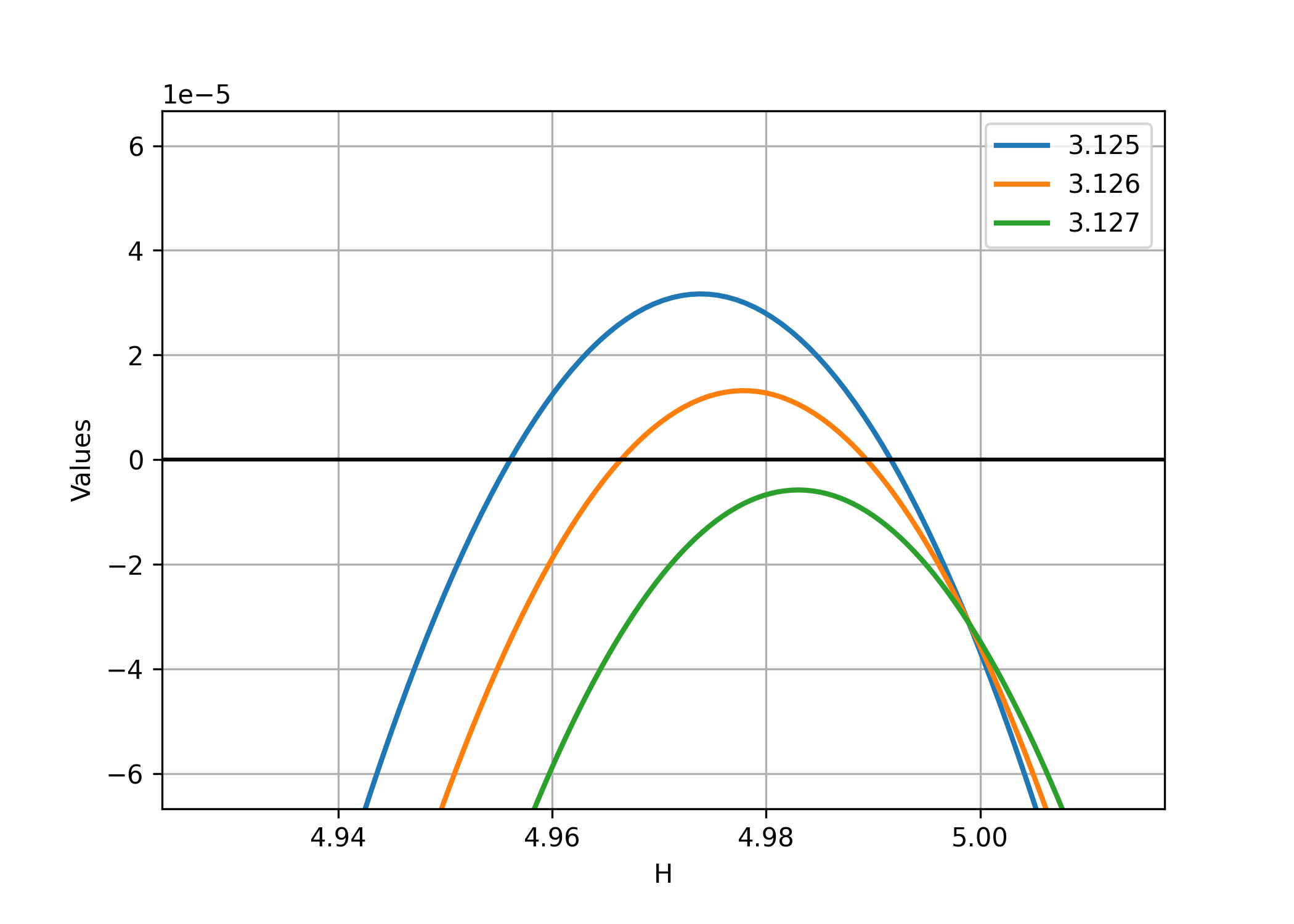} 
\includegraphics[width=6cm]{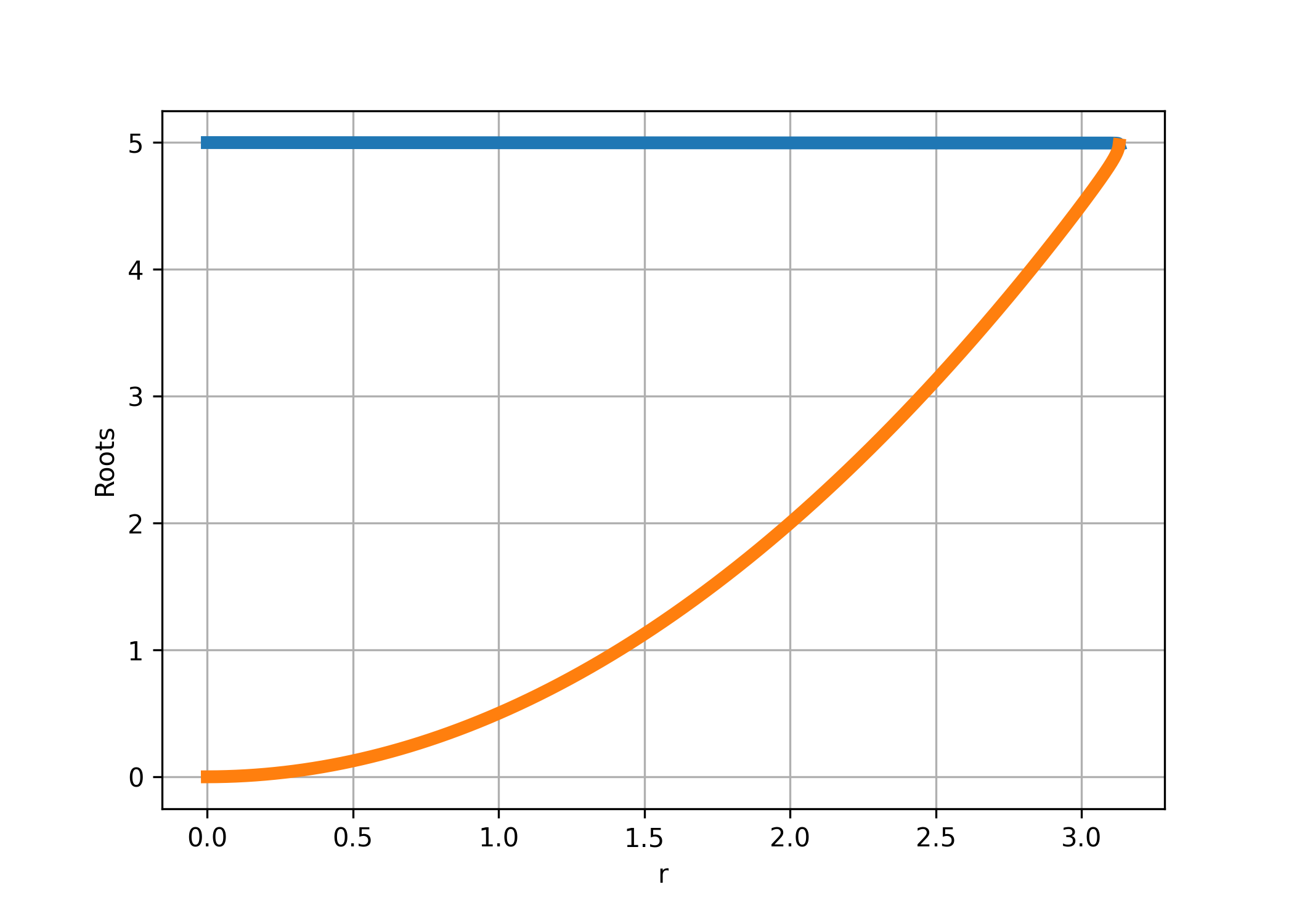} 
\includegraphics[width=6cm]{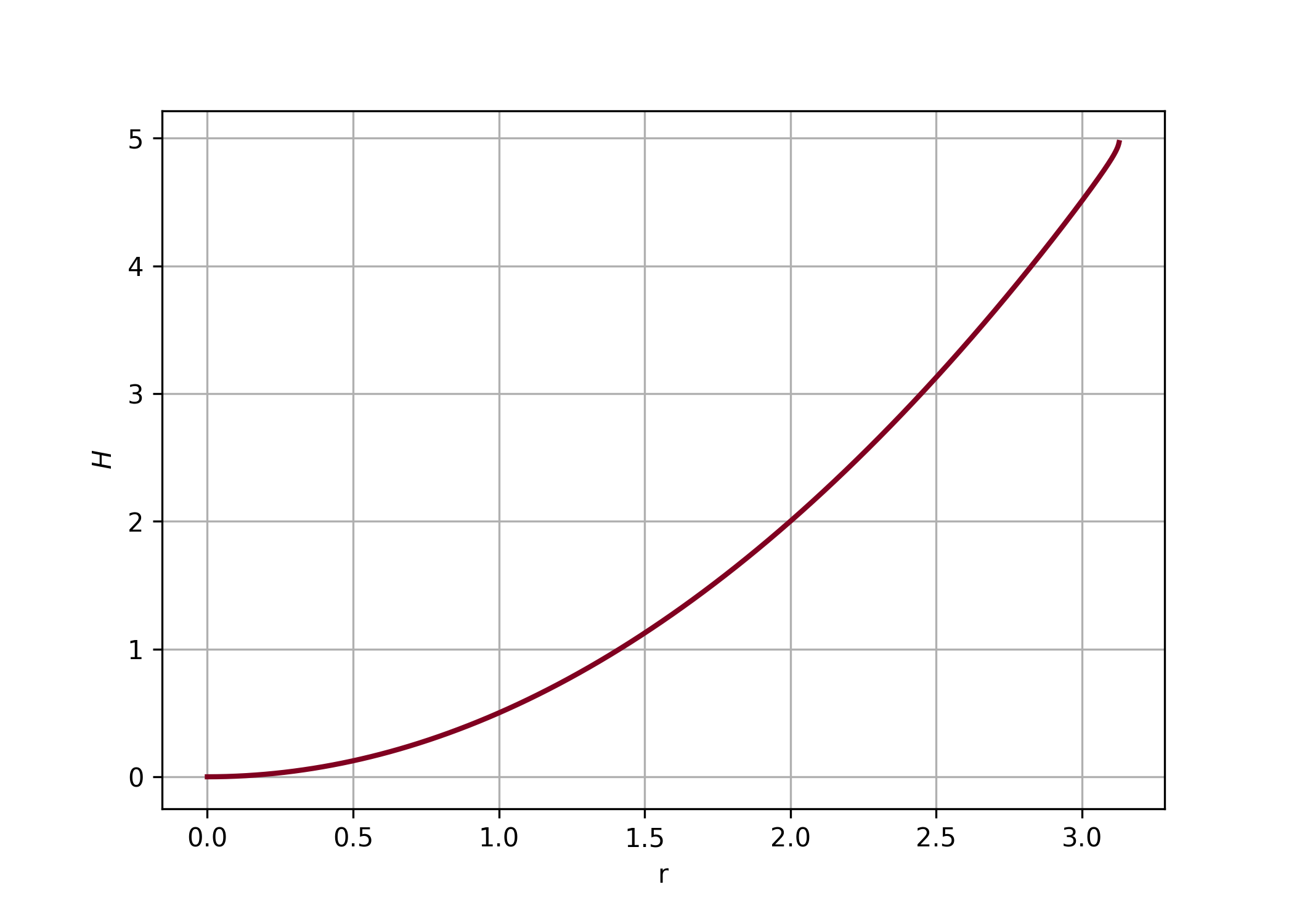} 
\includegraphics[width=6cm]{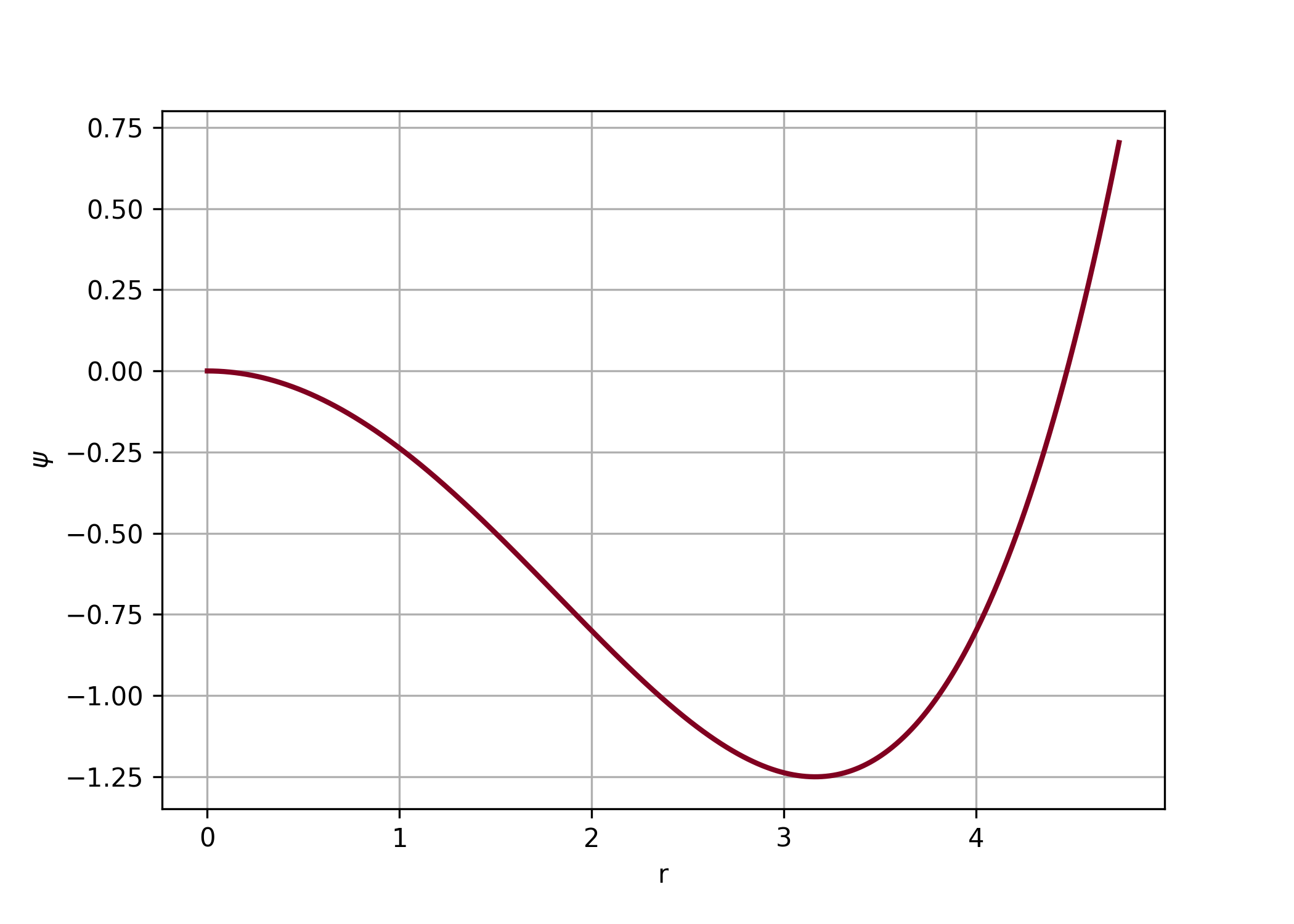} 
\caption{\textbf{Numerical integration of a radially symmetric system.}
(A) Solutions of Equation \ref{FD_Scheme_quad} for three consecutive steps ($\Delta r = 10^{-3}$) for Equation \ref{symm} (Parameter values: $\beta = 10$, $k = 0.05$, $\omega = -1$). Solutions vanish close to the limit cycle ($r = \sqrt{\beta} \approx 3.162$), as shown by the green parabola as does not cut the x-axis (black line).
(B) Evolution of the roots of Equation \ref{FD_Scheme_quad} as the system is integrated. 
(C) Solution for the Hamiltonian $H$ using Equation \ref{FD_Scheme_quad} up to the vanishing solution.
(D) Solution for the potential $\psi$, given by $p = 1$.}
\label{rad_fig}
\end{figure}

A first-order forward finite difference approximation of the derivative is given by
\begin{equation} \label{FDE}
\partial_r H \approx \frac{H_n - H_{n-1}}{\Delta r}
\end{equation}
By substituting Equation \ref{FDE} in Equation \ref{radial_polar}, a system of second degree algebraic equations is obtained, i.e.
\begin{equation} \label{FD_Scheme_quad}
-q H_n^2 +(p_0 + q H_{n-1}) H_n - (\dot r \Delta r + p_0 H_{n-1}) = 0,
\end{equation}
As the n-th value of the scalar function only depends on the previous one, the system can be solved by setting $H_0 = 0$ and solving for the following terms, one at a time. This is in fact the definition of the fixed point as the boundary condition and propagating the solution radially outwards.
Equation \ref{FD_Scheme_quad} has two solutions, which can be easily seen when solving the point $H_1$ (see Figure \ref{rad_fig}B), i.e.
\begin{equation} \label{FD_Scheme_1}
-q H_1^2 + p_0 H_1 - \dot r \Delta r = 0
\end{equation}
Equation \ref{FD_Scheme_quad} has real solutions when
\begin{equation} \label{FD_Scheme_sols}
(p_0 + q H_{n-1})^2\geq 4q(\dot r \Delta r + p_0 H_{n-1})
\end{equation}
If $H_1 = 1/2(\Delta r)^2$ and $\dot r = \Delta r(p_0 -q\Delta r)$ (as in Equation \ref{symm}), one can see that the solutions for Equation \ref{FD_Scheme_1} are $\Delta r$ and $p_0/q - \Delta r$, the latter being the solution for the limit cycle. When the solution approaches the limit cycle ($H_{n} \rightarrow \frac{p}{q}$, $p(H_n) \rightarrow 0$) this criteria is no longer met and no solution can be found (see Figure \ref{rad_fig}A and B).

Singularities might be a result of an implicit division by zero, when $p_0 - qH \approx 0$ and the solution blows up. Nonetheless, the parabolic structure of the solution is retrieved up to the limit cycle (see Figure \ref{rad_fig}C). As expected, the potential can be recovered by setting $p = 1$ (see Figure \ref{rad_fig}D).

\subsection{Asymmetric system}
In practice, systems do not exhibit the idealised properties of systems such as the one described by Equation \ref{symm}. Therefore, from now on, I will focus on solving the geometry of Equation \ref{asymm} using a simple forward method.

\begin{figure}[t]
\centering
\includegraphics[width=12cm]{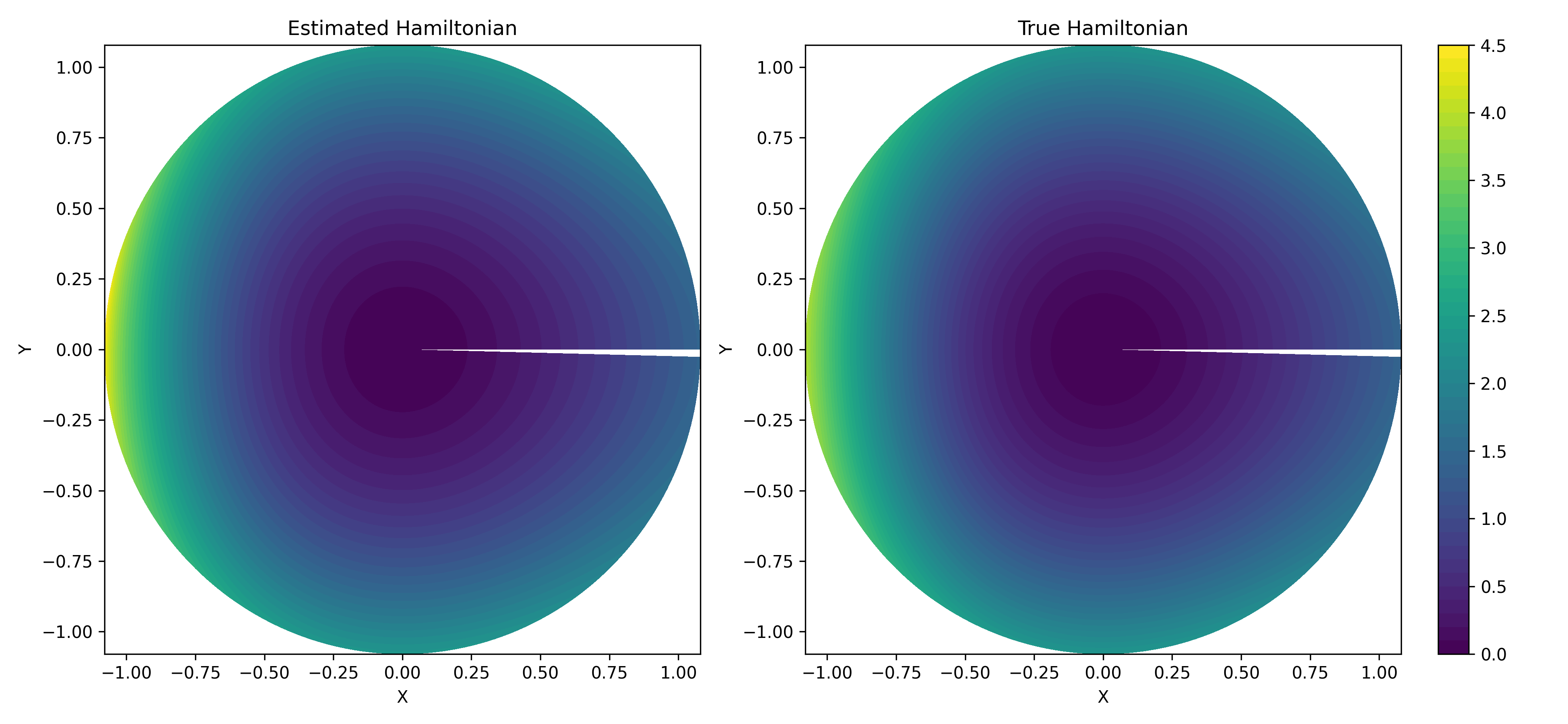} 
\includegraphics[width=6cm]{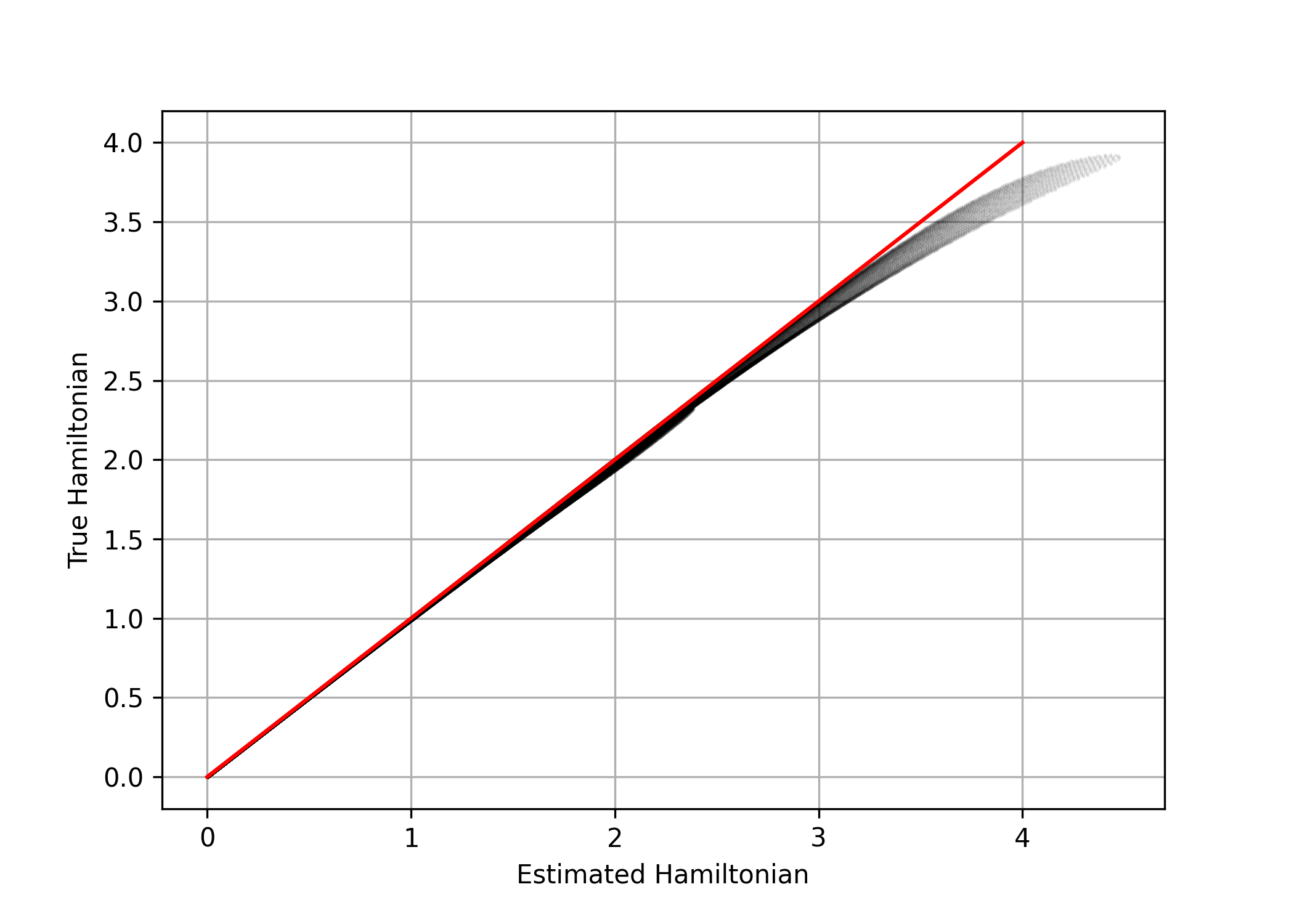} 
\caption{\textbf{Numerical integration of an asymmetric system.}
(A) Estimated and real Hamiltonians using the forward PDE method within the simulation horizon. Colorbar represents the value of Hamiltonian.
(B) Comparison between the real Hamiltonian and the estimation. Black dots represent each position and the red line represents perfect match. As the system gets farther away from the origin, the estimation gets worse.}
\label{pde_fig}
\end{figure}

As shown in Figure \ref{pde_fig}, the solution near the fixed point is close to the real value. However, as the system departs from the fixed point, the error increases. At a critical value, the method becomes unstable, and no solution is obtained. This is in line with the unstability of such forward methods.

Notwithstanding that this method is not suitable to get a complete solution, it can still be used to get a reasonable approximation of the dynamics near the fixed point and be the base of higher order methods. However, these methods might require additional conditions that one might not want to assume or impose.

\subsection{Computing potentials}
If $p = 1$, then Equation \ref{radial_polar} becomes
$
\dot r - \partial_r \psi = 0,
$
which can be solved as via Finite Differences (or in this case, the Euler Forward Method). As I defined $\psi = H(\beta - 2H)$, the Hamiltonian can be recovered as 
$
H = \frac{1}{2}(\beta \pm \sqrt{\beta^2 - 4 \psi}),
$
such the sign is chosen to make sure the correct form for $H$ is chosen.

\clearpage
\section{Analysis via Differential Geometry} \label{diff_geo}
In Einstein summation convention, Equation \ref{main_equation} becomes
\begin{equation} \label{main_ES}
\dot x^j = \kappa^{jk} \partial_k H,
\end{equation}
where $x^j$ is the $j$-th variable, and $\partial_k \equiv \partial /\partial x^k$.
In essence, Equation \ref{main_ES} reveals the invariance of the system when the generalised Hamiltonian $H$ or the bilinear from $\kappa$ is changed. As it has been discussed extensively in this paper, the choice of either one is arbitrary and the ideal one is based on the extraction of desired properties, such as the potential part vanishing at the limit cycle.

It is worth taking a look to special cases. If the dynamic tensor is diagonal, $\kappa^{jk} = g^{jk}$, then the system becomes a special case of the gradient descent on a (Riemannian) manifold \cite{rand2021}.
On the other hand, if $\kappa^{jk}$ is skew-symmetric, nondegenerate and with determinant one, we obtain the 2-dimensional symplectic form $\omega^{jk}$, which is the object of study of Hamiltonian systems using the language of differential geometry.
In summary, Equation \ref{main_ES} is some kind of generalisation of differential systems.

Interestingly, symplectic forms do not have any local invariant, while Riemannian manifolds do (the curvature). It would be interesting to see whether limit cycles, under the formulation of this paper, inherit local properties from the metric part, and thus, it could be possible to define limit cycles via a collection of local points with a certain invariant.

On this line, and using the language of differential forms, the problem here is finding a 2-covariant (normal) tensor that generates the exact form $dH$ when it is equiped to the differential system $\dot x^j$, i.e.
\begin{equation}
dH = \partial_k H d x^k = \kappa_{jk} \dot x^j d x^k,
\end{equation}
where $d$ is the exterior derivative of the 2-form $H$ and $\kappa_{jk}$ is the inverse of $\kappa^{jk}$. As $dH$ is exact, then $ddH = 0$ (for further read, see \cite{frankel2011}).

Although the exposition in this section is shallow, I hope this can spark the interest and further motivate research in this direction.

\end{document}